\documentclass[12pt]{amsart}
\title[classification of affine fully commutative elements]{A classification of affine fully commutative elements}
\author{Sadek AL HARBAT}
\address{IMJ Université Paris 7} 
\email{sadikharbat@math.univ-paris-diderot.fr}

\usepackage{etex} 
\usepackage[latin1]{inputenc}	
\usepackage[T1]{fontenc}
\usepackage[english]{babel}
\usepackage{amsmath,amssymb}
\usepackage{wasysym}
\usepackage{stmaryrd}
\usepackage[table]{xcolor}
\usepackage{color}
\usepackage{dsfont}\let\mathbb\mathds
\usepackage[all]{xy}
\usepackage{graphicx}
\usepackage{mathenv}
\usepackage{lastpage}
\usepackage{fancyhdr}	
\usepackage{subfigure}
\usepackage{geometry}
\usepackage[version=3]{mhchem}
\usepackage[force,almostfull]{textcomp}
\usepackage{array}
\usepackage{lipsum}		
\usepackage{eso-pic}
\usepackage{multirow}
\usepackage{fancybox}
\usepackage{cite}
\usepackage{amsthm}
\usepackage{listings}
\usepackage{lscape}
\usepackage{multicol}		
\usepackage{setspace}
\usepackage{times}
\usepackage{eurosym}
\usepackage{latexsym}
\usepackage{ulem}
\usepackage{lmodern}
\usepackage{colortbl}
\usepackage{rotating}
\usepackage{type1cm}
\usepackage{lettrine}
\usepackage{ragged2e}
\usepackage{mathrsfs}
\usepackage{scrtime}
\usepackage{enumitem}

\usepackage{longtable}
\usepackage{url}
\usepackage{nomentbl}
\usepackage{nomencl}
\usepackage{hyperref}
\hypersetup{
pdfpagemode=UseOutlines,      
pdfstartview=Fit,             
pdffitwindow=true,            
pdfpagelayout=TwoColumnsRight,
pdftoolbar=true,              
pdfmenubar=true,              
bookmarksopen=false,          
bookmarksnumbered=true,       
colorlinks=true,              
pdfauthor={Ton nom},          
pdftitle={Titre PDF},         
pdfcreator=PDFLaTeX,          %
pdfproducer=PDFLaTeX,         %
linkcolor=blue,               
urlcolor=blue,                
anchorcolor=black,            
citecolor=blue,               
frenchlinks=true,             
pdfborder={0 0 0}             
}
\usepackage{pgf, tikz, tikz-cd}
\usetikzlibrary{matrix,positioning}
\usetikzlibrary{arrows,decorations.markings}

\definecolor{aliceblue}{rgb}{0.94,0.97,1.00}
\definecolor{antiquewhite}{rgb}{0.98,0.92,0.84}
\definecolor{antiquewhite1}{rgb}{1.00,0.94,0.86}
\definecolor{antiquewhite2}{rgb}{0.93,0.87,0.80}
\definecolor{antiquewhite3}{rgb}{0.80,0.75,0.69}
\definecolor{antiquewhite4}{rgb}{0.55,0.51,0.47}
\definecolor{aquamarine}{rgb}{0.50,1.00,0.83}
\definecolor{aquamarine1}{rgb}{0.50,1.00,0.83}
\definecolor{aquamarine2}{rgb}{0.46,0.93,0.78}
\definecolor{aquamarine3}{rgb}{0.40,0.80,0.67}
\definecolor{aquamarine4}{rgb}{0.27,0.55,0.45}
\definecolor{azure}{rgb}{0.94,1.00,1.00}
\definecolor{azure1}{rgb}{0.94,1.00,1.00}
\definecolor{azure2}{rgb}{0.88,0.93,0.93}
\definecolor{azure3}{rgb}{0.76,0.80,0.80}
\definecolor{azure4}{rgb}{0.51,0.55,0.55}
\definecolor{beige}{rgb}{0.96,0.96,0.86}
\definecolor{bisque}{rgb}{1.00,0.89,0.77}
\definecolor{bisque1}{rgb}{1.00,0.89,0.77}
\definecolor{bisque2}{rgb}{0.93,0.84,0.72}
\definecolor{bisque3}{rgb}{0.80,0.72,0.62}
\definecolor{bisque4}{rgb}{0.55,0.49,0.42}
\definecolor{black}{rgb}{0.00,0.00,0.00}
\definecolor{blanchedalmond}{rgb}{1.00,0.92,0.80}
\definecolor{blue}{rgb}{0.00,0.00,1.00}
\definecolor{blue1}{rgb}{0.00,0.00,1.00}
\definecolor{blue2}{rgb}{0.00,0.00,0.93}
\definecolor{blue3}{rgb}{0.00,0.00,0.80}
\definecolor{blue4}{rgb}{0.00,0.00,0.55}
\definecolor{blueviolet}{rgb}{0.54,0.17,0.89}
\definecolor{brown}{rgb}{0.65,0.16,0.16}
\definecolor{brown1}{rgb}{1.00,0.25,0.25}
\definecolor{brown2}{rgb}{0.93,0.23,0.23}
\definecolor{brown3}{rgb}{0.80,0.20,0.20}
\definecolor{brown4}{rgb}{0.55,0.14,0.14}
\definecolor{burlywood}{rgb}{0.87,0.72,0.53}
\definecolor{burlywood1}{rgb}{1.00,0.83,0.61}
\definecolor{burlywood2}{rgb}{0.93,0.77,0.57}
\definecolor{burlywood3}{rgb}{0.80,0.67,0.49}
\definecolor{burlywood4}{rgb}{0.55,0.45,0.33}
\definecolor{cadetblue}{rgb}{0.37,0.62,0.63}
\definecolor{cadetblue1}{rgb}{0.60,0.96,1.00}
\definecolor{cadetblue2}{rgb}{0.56,0.90,0.93}
\definecolor{cadetblue3}{rgb}{0.48,0.77,0.80}
\definecolor{cadetblue4}{rgb}{0.33,0.53,0.55}
\definecolor{chartreuse}{rgb}{0.50,1.00,0.00}
\definecolor{chartreuse1}{rgb}{0.50,1.00,0.00}
\definecolor{chartreuse2}{rgb}{0.46,0.93,0.00}
\definecolor{chartreuse3}{rgb}{0.40,0.80,0.00}
\definecolor{chartreuse4}{rgb}{0.27,0.55,0.00}
\definecolor{chocolate}{rgb}{0.82,0.41,0.12}
\definecolor{chocolate1}{rgb}{1.00,0.50,0.14}
\definecolor{chocolate2}{rgb}{0.93,0.46,0.13}
\definecolor{chocolate3}{rgb}{0.80,0.40,0.11}
\definecolor{chocolate4}{rgb}{0.55,0.27,0.07}
\definecolor{coral}{rgb}{1.00,0.50,0.31}
\definecolor{coral1}{rgb}{1.00,0.45,0.34}
\definecolor{coral2}{rgb}{0.93,0.42,0.31}
\definecolor{coral3}{rgb}{0.80,0.36,0.27}
\definecolor{coral4}{rgb}{0.55,0.24,0.18}
\definecolor{cornsilk}{rgb}{1.00,0.97,0.86}
\definecolor{cornflowerblue}{rgb}{0.39,0.58,0.93}
\definecolor{cornsilk1}{rgb}{1.00,0.97,0.86}
\definecolor{cornsilk2}{rgb}{0.93,0.91,0.80}
\definecolor{cornsilk3}{rgb}{0.80,0.78,0.69}
\definecolor{cornsilk4}{rgb}{0.55,0.53,0.47}
\definecolor{cyan}{rgb}{0.00,1.00,1.00}
\definecolor{cyan1}{rgb}{0.00,1.00,1.00}
\definecolor{cyan2}{rgb}{0.00,0.93,0.93}
\definecolor{cyan3}{rgb}{0.00,0.80,0.80}
\definecolor{cyan4}{rgb}{0.00,0.55,0.55}
\definecolor{darkblue}{rgb}{0.00,0.00,0.55}
\definecolor{darkcyan}{rgb}{0.00,0.55,0.55}
\definecolor{darkgoldenrod}{rgb}{0.72,0.53,0.04}
\definecolor{DarkGoldenrod1}{rgb}{1.00,0.73,0.06}
\definecolor{DarkGoldenrod2}{rgb}{0.93,0.68,0.05}
\definecolor{DarkGoldenrod3}{rgb}{0.80,0.58,0.05}
\definecolor{DarkGoldenrod4}{rgb}{0.55,0.40,0.03}
\definecolor{darkgray}{rgb}{0.66,0.66,0.66}
\definecolor{darkgreen}{rgb}{0.00,0.39,0.00}
\definecolor{darkgrey}{rgb}{0.66,0.66,0.66}
\definecolor{darkkhaki}{rgb}{0.74,0.72,0.42}
\definecolor{darkmagenta}{rgb}{0.55,0.00,0.55}
\definecolor{darkOlivegreen}{rgb}{0.33,0.42,0.18}
\definecolor{darkOlivegreen1}{rgb}{0.79,1.00,0.44}
\definecolor{darkOlivegreen2}{rgb}{0.74,0.93,0.41}
\definecolor{darkOlivegreen3}{rgb}{0.64,0.80,0.35}
\definecolor{darkOlivegreen4}{rgb}{0.43,0.55,0.24}
\definecolor{darkolive}{rgb}{0.33,0.42,0.18}
\definecolor{DarkOrange1}{rgb}{1.00,0.50,0.00}
\definecolor{DarkOrange2}{rgb}{0.93,0.46,0.00}
\definecolor{DarkOrange3}{rgb}{0.80,0.40,0.00}
\definecolor{DarkOrange4}{rgb}{0.55,0.27,0.00}
\definecolor{darkorange}{rgb}{1.00,0.55,0.00}
\definecolor{darkorchid}{rgb}{0.60,0.20,0.80}
\definecolor{darkOrchid1}{rgb}{0.75,0.24,1.00}
\definecolor{darkOrchid2}{rgb}{0.70,0.23,0.93}
\definecolor{darkOrchid3}{rgb}{0.60,0.20,0.80}
\definecolor{darkOrchid4}{rgb}{0.41,0.13,0.55}
\definecolor{darkred}{rgb}{0.55,0.00,0.00}
\definecolor{darksalmon}{rgb}{0.91,0.59,0.48}
\definecolor{darkseagreen}{rgb}{0.56,0.74,0.56}
\definecolor{darkseagreen1}{rgb}{0.76,1.00,0.76}
\definecolor{darkseagreen2}{rgb}{0.71,0.93,0.71}
\definecolor{darkseagreen3}{rgb}{0.61,0.80,0.61}
\definecolor{darkseagreen4}{rgb}{0.41,0.55,0.41}
\definecolor{darksea}{rgb}{0.56,0.74,0.56}
\definecolor{darkslategray}{rgb}{0.18,0.31,0.31}
\definecolor{darkslateblue}{rgb}{0.28,0.24,0.55}
\definecolor{darkslategray1}{rgb}{0.59,1.00,1.00}
\definecolor{darkslategray2}{rgb}{0.55,0.93,0.93}
\definecolor{darkslategray3}{rgb}{0.47,0.80,0.80}
\definecolor{darkslategray4}{rgb}{0.32,0.55,0.55}
\definecolor{darkslate}{rgb}{0.18,0.31,0.31}
\definecolor{darkslate1}{rgb}{0.28,0.24,0.55}
\definecolor{darkturquoise}{rgb}{0.00,0.81,0.82}
\definecolor{darkviolet}{rgb}{0.58,0.00,0.83}
\definecolor{deeppink}{rgb}{1.00,0.08,0.58}
\definecolor{deepPink1}{rgb}{1.00,0.08,0.58}
\definecolor{deepPink2}{rgb}{0.93,0.07,0.54}
\definecolor{deepPink3}{rgb}{0.80,0.06,0.46}
\definecolor{deepPink4}{rgb}{0.55,0.04,0.31}
\definecolor{deepskyblue}{rgb}{0.00,0.75,1.00}
\definecolor{deepskyblue1}{rgb}{0.00,0.75,1.00}
\definecolor{deepskyblue2}{rgb}{0.00,0.70,0.93}
\definecolor{deepskyblue3}{rgb}{0.00,0.60,0.80}
\definecolor{deepskyblue4}{rgb}{0.00,0.41,0.55}
\definecolor{deepsky}{rgb}{0.00,0.75,1.00}
\definecolor{dimgray}{rgb}{0.41,0.41,0.41}
\definecolor{dodgerblue}{rgb}{0.12,0.56,1.00}
\definecolor{dodgerblue1}{rgb}{0.12,0.56,1.00}
\definecolor{dodgerblue2}{rgb}{0.11,0.53,0.93}
\definecolor{dodgerblue3}{rgb}{0.09,0.45,0.80}
\definecolor{dodgerblue4}{rgb}{0.06,0.31,0.55}
\definecolor{firebrick}{rgb}{0.70,0.13,0.13}
\definecolor{firebrick1}{rgb}{1.00,0.19,0.19}
\definecolor{firebrick2}{rgb}{0.93,0.17,0.17}
\definecolor{firebrick3}{rgb}{0.80,0.15,0.15}
\definecolor{firebrick4}{rgb}{0.55,0.10,0.10}
\definecolor{floralwhite}{rgb}{1.00,0.98,0.94}
\definecolor{forestgreen}{rgb}{0.13,0.55,0.13}
\definecolor{gainsboro}{rgb}{0.86,0.86,0.86}
\definecolor{ghostwhite}{rgb}{0.97,0.97,1.00}
\definecolor{gold1}{rgb}{1.00,0.84,0.00}
\definecolor{gold2}{rgb}{0.93,0.79,0.00}
\definecolor{gold3}{rgb}{0.80,0.68,0.00}
\definecolor{gold4}{rgb}{0.55,0.46,0.00}
\definecolor{goldenrod}{rgb}{0.85,0.65,0.13}
\definecolor{goldenrod1}{rgb}{1.00,0.76,0.15}
\definecolor{goldenrod2}{rgb}{0.93,0.71,0.13}
\definecolor{goldenrod3}{rgb}{0.80,0.61,0.11}
\definecolor{goldenrod4}{rgb}{0.55,0.41,0.08}
\definecolor{gold}{rgb}{1.00,0.84,0.00}
\definecolor{gray}{rgb}{0.75,0.75,0.75}
\definecolor{gray0}{rgb}{0.00,0.00,0.00}
\definecolor{gray1}{rgb}{0.01,0.01,0.01}
\definecolor{gray2}{rgb}{0.02,0.02,0.02}
\definecolor{gray3}{rgb}{0.03,0.03,0.03}
\definecolor{gray4}{rgb}{0.04,0.04,0.04}
\definecolor{gray5}{rgb}{0.05,0.05,0.05}
\definecolor{gray6}{rgb}{0.06,0.06,0.06}
\definecolor{gray7}{rgb}{0.07,0.07,0.07}
\definecolor{gray8}{rgb}{0.08,0.08,0.08}
\definecolor{gray9}{rgb}{0.09,0.09,0.09}
\definecolor{gray10}{rgb}{0.10,0.10,0.10}
\definecolor{gray11}{rgb}{0.11,0.11,0.11}
\definecolor{gray12}{rgb}{0.12,0.12,0.12}
\definecolor{gray13}{rgb}{0.13,0.13,0.13}
\definecolor{gray14}{rgb}{0.14,0.14,0.14}
\definecolor{gray15}{rgb}{0.15,0.15,0.15}
\definecolor{gray16}{rgb}{0.16,0.16,0.16}
\definecolor{gray17}{rgb}{0.17,0.17,0.17}
\definecolor{gray18}{rgb}{0.18,0.18,0.18}
\definecolor{gray19}{rgb}{0.19,0.19,0.19}
\definecolor{gray20}{rgb}{0.20,0.20,0.20}
\definecolor{gray21}{rgb}{0.21,0.21,0.21}
\definecolor{gray22}{rgb}{0.22,0.22,0.22}
\definecolor{gray23}{rgb}{0.23,0.23,0.23}
\definecolor{gray24}{rgb}{0.24,0.24,0.24}
\definecolor{gray25}{rgb}{0.25,0.25,0.25}
\definecolor{gray26}{rgb}{0.26,0.26,0.26}
\definecolor{gray27}{rgb}{0.27,0.27,0.27}
\definecolor{gray28}{rgb}{0.28,0.28,0.28}
\definecolor{gray29}{rgb}{0.29,0.29,0.29}
\definecolor{gray30}{rgb}{0.30,0.30,0.30}
\definecolor{gray31}{rgb}{0.31,0.31,0.31}
\definecolor{gray32}{rgb}{0.32,0.32,0.32}
\definecolor{gray33}{rgb}{0.33,0.33,0.33}
\definecolor{gray34}{rgb}{0.34,0.34,0.34}
\definecolor{gray35}{rgb}{0.35,0.35,0.35}
\definecolor{gray36}{rgb}{0.36,0.36,0.36}
\definecolor{gray37}{rgb}{0.37,0.37,0.37}
\definecolor{gray38}{rgb}{0.38,0.38,0.38}
\definecolor{gray39}{rgb}{0.39,0.39,0.39}
\definecolor{gray40}{rgb}{0.40,0.40,0.40}
\definecolor{gray41}{rgb}{0.41,0.41,0.41}
\definecolor{gray42}{rgb}{0.42,0.42,0.42}
\definecolor{gray43}{rgb}{0.43,0.43,0.43}
\definecolor{gray44}{rgb}{0.44,0.44,0.44}
\definecolor{gray45}{rgb}{0.45,0.45,0.45}
\definecolor{gray46}{rgb}{0.46,0.46,0.46}
\definecolor{gray47}{rgb}{0.47,0.47,0.47}
\definecolor{gray48}{rgb}{0.48,0.48,0.48}
\definecolor{gray49}{rgb}{0.49,0.49,0.49}
\definecolor{gray50}{rgb}{0.50,0.50,0.50}
\definecolor{gray51}{rgb}{0.51,0.51,0.51}
\definecolor{gray52}{rgb}{0.52,0.52,0.52}
\definecolor{gray53}{rgb}{0.53,0.53,0.53}
\definecolor{gray54}{rgb}{0.54,0.54,0.54}
\definecolor{gray55}{rgb}{0.55,0.55,0.55}
\definecolor{gray56}{rgb}{0.56,0.56,0.56}
\definecolor{gray57}{rgb}{0.57,0.57,0.57}
\definecolor{gray58}{rgb}{0.58,0.58,0.58}
\definecolor{gray59}{rgb}{0.59,0.59,0.59}
\definecolor{gray60}{rgb}{0.60,0.60,0.60}
\definecolor{gray61}{rgb}{0.61,0.61,0.61}
\definecolor{gray62}{rgb}{0.62,0.62,0.62}
\definecolor{gray63}{rgb}{0.63,0.63,0.63}
\definecolor{gray64}{rgb}{0.64,0.64,0.64}
\definecolor{gray65}{rgb}{0.65,0.65,0.65}
\definecolor{gray66}{rgb}{0.66,0.66,0.66}
\definecolor{gray67}{rgb}{0.67,0.67,0.67}
\definecolor{gray68}{rgb}{0.68,0.68,0.68}
\definecolor{gray69}{rgb}{0.69,0.69,0.69}
\definecolor{gray70}{rgb}{0.70,0.70,0.70}
\definecolor{gray71}{rgb}{0.71,0.71,0.71}
\definecolor{gray72}{rgb}{0.72,0.72,0.72}
\definecolor{gray73}{rgb}{0.73,0.73,0.73}
\definecolor{gray74}{rgb}{0.74,0.74,0.74}
\definecolor{gray75}{rgb}{0.75,0.75,0.75}
\definecolor{gray76}{rgb}{0.76,0.76,0.76}
\definecolor{gray77}{rgb}{0.77,0.77,0.77}
\definecolor{gray78}{rgb}{0.78,0.78,0.78}
\definecolor{gray79}{rgb}{0.79,0.79,0.79}
\definecolor{gray80}{rgb}{0.80,0.80,0.80}
\definecolor{gray81}{rgb}{0.81,0.81,0.81}
\definecolor{gray82}{rgb}{0.82,0.82,0.82}
\definecolor{gray83}{rgb}{0.83,0.83,0.83}
\definecolor{gray84}{rgb}{0.84,0.84,0.84}
\definecolor{gray85}{rgb}{0.85,0.85,0.85}
\definecolor{gray86}{rgb}{0.86,0.86,0.86}
\definecolor{gray87}{rgb}{0.87,0.87,0.87}
\definecolor{gray88}{rgb}{0.88,0.88,0.88}
\definecolor{gray89}{rgb}{0.89,0.89,0.89}
\definecolor{gray90}{rgb}{0.90,0.90,0.90}
\definecolor{gray91}{rgb}{0.91,0.91,0.91}
\definecolor{gray92}{rgb}{0.92,0.92,0.92}
\definecolor{gray93}{rgb}{0.93,0.93,0.93}
\definecolor{gray94}{rgb}{0.94,0.94,0.94}
\definecolor{gray95}{rgb}{0.95,0.95,0.95}
\definecolor{gray96}{rgb}{0.96,0.96,0.96}
\definecolor{gray97}{rgb}{0.97,0.97,0.97}
\definecolor{gray98}{rgb}{0.98,0.98,0.98}
\definecolor{gray99}{rgb}{0.99,0.99,0.99}
\definecolor{gray100}{rgb}{1.00,1.00,1.00}
\definecolor{green}{rgb}{0.00,1.00,0.00}
\definecolor{green1}{rgb}{0.00,1.00,0.00}
\definecolor{green2}{rgb}{0.00,0.93,0.00}
\definecolor{green3}{rgb}{0.00,0.80,0.00}
\definecolor{green4}{rgb}{0.00,0.55,0.00}
\definecolor{greenyellow}{rgb}{0.68,1.00,0.18}
\definecolor{grey}{rgb}{0.75,0.75,0.75}
\definecolor{grey0}{rgb}{0.00,0.00,0.00}
\definecolor{grey1}{rgb}{0.01,0.01,0.01}
\definecolor{grey2}{rgb}{0.02,0.02,0.02}
\definecolor{grey3}{rgb}{0.03,0.03,0.03}
\definecolor{grey4}{rgb}{0.04,0.04,0.04}
\definecolor{grey5}{rgb}{0.05,0.05,0.05}
\definecolor{grey6}{rgb}{0.06,0.06,0.06}
\definecolor{grey7}{rgb}{0.07,0.07,0.07}
\definecolor{grey8}{rgb}{0.08,0.08,0.08}
\definecolor{grey9}{rgb}{0.09,0.09,0.09}
\definecolor{grey10}{rgb}{0.10,0.10,0.10}
\definecolor{grey11}{rgb}{0.11,0.11,0.11}
\definecolor{grey12}{rgb}{0.12,0.12,0.12}
\definecolor{grey13}{rgb}{0.13,0.13,0.13}
\definecolor{grey14}{rgb}{0.14,0.14,0.14}
\definecolor{grey15}{rgb}{0.15,0.15,0.15}
\definecolor{grey16}{rgb}{0.16,0.16,0.16}
\definecolor{grey17}{rgb}{0.17,0.17,0.17}
\definecolor{grey18}{rgb}{0.18,0.18,0.18}
\definecolor{grey19}{rgb}{0.19,0.19,0.19}
\definecolor{grey20}{rgb}{0.20,0.20,0.20}
\definecolor{grey21}{rgb}{0.21,0.21,0.21}
\definecolor{grey22}{rgb}{0.22,0.22,0.22}
\definecolor{grey23}{rgb}{0.23,0.23,0.23}
\definecolor{grey24}{rgb}{0.24,0.24,0.24}
\definecolor{grey25}{rgb}{0.25,0.25,0.25}
\definecolor{grey26}{rgb}{0.26,0.26,0.26}
\definecolor{grey27}{rgb}{0.27,0.27,0.27}
\definecolor{grey28}{rgb}{0.28,0.28,0.28}
\definecolor{grey29}{rgb}{0.29,0.29,0.29}
\definecolor{grey30}{rgb}{0.30,0.30,0.30}
\definecolor{grey31}{rgb}{0.31,0.31,0.31}
\definecolor{grey32}{rgb}{0.32,0.32,0.32}
\definecolor{grey33}{rgb}{0.33,0.33,0.33}
\definecolor{grey34}{rgb}{0.34,0.34,0.34}
\definecolor{grey35}{rgb}{0.35,0.35,0.35}
\definecolor{grey36}{rgb}{0.36,0.36,0.36}
\definecolor{grey37}{rgb}{0.37,0.37,0.37}
\definecolor{grey38}{rgb}{0.38,0.38,0.38}
\definecolor{grey39}{rgb}{0.39,0.39,0.39}
\definecolor{grey40}{rgb}{0.40,0.40,0.40}
\definecolor{grey41}{rgb}{0.41,0.41,0.41}
\definecolor{grey42}{rgb}{0.42,0.42,0.42}
\definecolor{grey43}{rgb}{0.43,0.43,0.43}
\definecolor{grey44}{rgb}{0.44,0.44,0.44}
\definecolor{grey45}{rgb}{0.45,0.45,0.45}
\definecolor{grey46}{rgb}{0.46,0.46,0.46}
\definecolor{grey47}{rgb}{0.47,0.47,0.47}
\definecolor{grey48}{rgb}{0.48,0.48,0.48}
\definecolor{grey49}{rgb}{0.49,0.49,0.49}
\definecolor{grey50}{rgb}{0.50,0.50,0.50}
\definecolor{grey51}{rgb}{0.51,0.51,0.51}
\definecolor{grey52}{rgb}{0.52,0.52,0.52}
\definecolor{grey53}{rgb}{0.53,0.53,0.53}
\definecolor{grey54}{rgb}{0.54,0.54,0.54}
\definecolor{grey55}{rgb}{0.55,0.55,0.55}
\definecolor{grey56}{rgb}{0.56,0.56,0.56}
\definecolor{grey57}{rgb}{0.57,0.57,0.57}
\definecolor{grey58}{rgb}{0.58,0.58,0.58}
\definecolor{grey59}{rgb}{0.59,0.59,0.59}
\definecolor{grey60}{rgb}{0.60,0.60,0.60}
\definecolor{grey61}{rgb}{0.61,0.61,0.61}
\definecolor{grey62}{rgb}{0.62,0.62,0.62}
\definecolor{grey63}{rgb}{0.63,0.63,0.63}
\definecolor{grey64}{rgb}{0.64,0.64,0.64}
\definecolor{grey65}{rgb}{0.65,0.65,0.65}
\definecolor{grey66}{rgb}{0.66,0.66,0.66}
\definecolor{grey67}{rgb}{0.67,0.67,0.67}
\definecolor{grey68}{rgb}{0.68,0.68,0.68}
\definecolor{grey69}{rgb}{0.69,0.69,0.69}
\definecolor{grey70}{rgb}{0.70,0.70,0.70}
\definecolor{grey71}{rgb}{0.71,0.71,0.71}
\definecolor{grey72}{rgb}{0.72,0.72,0.72}
\definecolor{grey73}{rgb}{0.73,0.73,0.73}
\definecolor{grey74}{rgb}{0.74,0.74,0.74}
\definecolor{grey75}{rgb}{0.75,0.75,0.75}
\definecolor{grey76}{rgb}{0.76,0.76,0.76}
\definecolor{grey77}{rgb}{0.77,0.77,0.77}
\definecolor{grey78}{rgb}{0.78,0.78,0.78}
\definecolor{grey79}{rgb}{0.79,0.79,0.79}
\definecolor{grey80}{rgb}{0.80,0.80,0.80}
\definecolor{grey81}{rgb}{0.81,0.81,0.81}
\definecolor{grey82}{rgb}{0.82,0.82,0.82}
\definecolor{grey83}{rgb}{0.83,0.83,0.83}
\definecolor{grey84}{rgb}{0.84,0.84,0.84}
\definecolor{grey85}{rgb}{0.85,0.85,0.85}
\definecolor{grey86}{rgb}{0.86,0.86,0.86}
\definecolor{grey87}{rgb}{0.87,0.87,0.87}
\definecolor{grey88}{rgb}{0.88,0.88,0.88}
\definecolor{grey89}{rgb}{0.89,0.89,0.89}
\definecolor{grey90}{rgb}{0.90,0.90,0.90}
\definecolor{grey91}{rgb}{0.91,0.91,0.91}
\definecolor{grey92}{rgb}{0.92,0.92,0.92}
\definecolor{grey93}{rgb}{0.93,0.93,0.93}
\definecolor{grey94}{rgb}{0.94,0.94,0.94}
\definecolor{grey95}{rgb}{0.95,0.95,0.95}
\definecolor{grey96}{rgb}{0.96,0.96,0.96}
\definecolor{grey97}{rgb}{0.97,0.97,0.97}
\definecolor{grey98}{rgb}{0.98,0.98,0.98}
\definecolor{grey99}{rgb}{0.99,0.99,0.99}
\definecolor{grey100}{rgb}{1.00,1.00,1.00}
\definecolor{honeydew}{rgb}{0.94,1.00,0.94}
\definecolor{honeydew1}{rgb}{0.94,1.00,0.94}
\definecolor{honeydew2}{rgb}{0.88,0.93,0.88}
\definecolor{honeydew3}{rgb}{0.76,0.80,0.76}
\definecolor{honeydew4}{rgb}{0.51,0.55,0.51}
\definecolor{hotpink}{rgb}{1.00,0.41,0.71}
\definecolor{hotPink1}{rgb}{1.00,0.43,0.71}
\definecolor{hotPink2}{rgb}{0.93,0.42,0.65}
\definecolor{hotPink3}{rgb}{0.80,0.38,0.56}
\definecolor{hotPink4}{rgb}{0.55,0.23,0.38}
\definecolor{indianred}{rgb}{0.80,0.36,0.36}
\definecolor{indianred1}{rgb}{1.00,0.42,0.42}
\definecolor{indianred2}{rgb}{0.93,0.39,0.39}
\definecolor{indianred3}{rgb}{0.80,0.33,0.33}
\definecolor{indianred4}{rgb}{0.55,0.23,0.23}
\definecolor{ivory}{rgb}{1.00,1.00,0.94}
\definecolor{ivory1}{rgb}{1.00,1.00,0.94}
\definecolor{ivory2}{rgb}{0.93,0.93,0.88}
\definecolor{ivory3}{rgb}{0.80,0.80,0.76}
\definecolor{ivory4}{rgb}{0.55,0.55,0.51}
\definecolor{khaki}{rgb}{0.94,0.90,0.55}
\definecolor{khaki1}{rgb}{1.00,0.96,0.56}
\definecolor{khaki2}{rgb}{0.93,0.90,0.52}
\definecolor{khaki3}{rgb}{0.80,0.78,0.45}
\definecolor{khaki4}{rgb}{0.55,0.53,0.31}
\definecolor{lavenderblush}{rgb}{1.00,0.94,0.96}
\definecolor{lavenderblush1}{rgb}{1.00,0.94,0.96}
\definecolor{lavenderblush2}{rgb}{0.93,0.88,0.90}
\definecolor{lavenderblush3}{rgb}{0.80,0.76,0.77}
\definecolor{lavenderblush4}{rgb}{0.55,0.51,0.53}
\definecolor{lavender}{rgb}{0.90,0.90,0.98}
\definecolor{lawngreen}{rgb}{0.49,0.99,0.00}
\definecolor{lemonchiffon}{rgb}{1.00,0.98,0.80}
\definecolor{lemonchiffon1}{rgb}{1.00,0.98,0.80}
\definecolor{lemonchiffon2}{rgb}{0.93,0.91,0.75}
\definecolor{lemonchiffon3}{rgb}{0.80,0.79,0.65}
\definecolor{lemonchiffon4}{rgb}{0.55,0.54,0.44}
\definecolor{lightblue}{rgb}{0.68,0.85,0.90}
\definecolor{lightblue1}{rgb}{0.75,0.94,1.00}
\definecolor{lightblue2}{rgb}{0.70,0.87,0.93}
\definecolor{lightblue3}{rgb}{0.60,0.75,0.80}
\definecolor{lightblue4}{rgb}{0.41,0.51,0.55}
\definecolor{lightcoral}{rgb}{0.94,0.50,0.50}
\definecolor{lightcyan}{rgb}{0.88,1.00,1.00}
\definecolor{lightcyan1}{rgb}{0.88,1.00,1.00}
\definecolor{lightcyan2}{rgb}{0.82,0.93,0.93}
\definecolor{lightcyan3}{rgb}{0.71,0.80,0.80}
\definecolor{lightcyan4}{rgb}{0.48,0.55,0.55}
\definecolor{lightgoldenrod}{rgb}{0.93,0.87,0.51}
\definecolor{lightgoldenrod0}{rgb}{0.98,0.98,0.82}
\definecolor{lightgoldenrod1}{rgb}{1.00,0.93,0.55}
\definecolor{lightgoldenrod2}{rgb}{0.93,0.86,0.51}
\definecolor{lightgoldenrod3}{rgb}{0.80,0.75,0.44}
\definecolor{lightgoldenrod4}{rgb}{0.55,0.51,0.30}
\definecolor{lightgoldenrodYellow}{rgb}{0.98,0.98,0.82}
\definecolor{lightgray}{rgb}{0.83,0.83,0.83}
\definecolor{lightgreen}{rgb}{0.56,0.93,0.56}
\definecolor{lightgrey}{rgb}{0.83,0.83,0.83}
\definecolor{lightpink}{rgb}{1.00,0.71,0.76}
\definecolor{lightpink1}{rgb}{1.00,0.68,0.73}
\definecolor{lightpink2}{rgb}{0.93,0.64,0.68}
\definecolor{lightpink3}{rgb}{0.80,0.55,0.58}
\definecolor{lightpink4}{rgb}{0.55,0.37,0.40}
\definecolor{lightsalmon}{rgb}{1.00,0.63,0.48}
\definecolor{lightsalmon1}{rgb}{1.00,0.63,0.48}
\definecolor{lightsalmon2}{rgb}{0.93,0.58,0.45}
\definecolor{lightsalmon3}{rgb}{0.80,0.51,0.38}
\definecolor{lightsalmon4}{rgb}{0.55,0.34,0.26}
\definecolor{lightseagreen}{rgb}{0.13,0.70,0.67}
\definecolor{lightsea}{rgb}{0.13,0.70,0.67}
\definecolor{lightsky}{rgb}{0.53,0.81,0.98}
\definecolor{lightSkyblue}{rgb}{0.53,0.81,0.98}
\definecolor{lightSkyblue1}{rgb}{0.69,0.89,1.00}
\definecolor{lightSkyblue2}{rgb}{0.64,0.83,0.93}
\definecolor{lightSkyblue3}{rgb}{0.55,0.71,0.80}
\definecolor{lightSkyblue4}{rgb}{0.38,0.48,0.55}
\definecolor{lightslateblue}{rgb}{0.52,0.44,1.00}
\definecolor{lightslategray}{rgb}{0.47,0.53,0.60}
\definecolor{lightslate}{rgb}{0.47,0.53,0.60}
\definecolor{lightslate1}{rgb}{0.52,0.44,1.00}
\definecolor{lightsteelblue}{rgb}{0.69,0.77,0.87}
\definecolor{lightsteelblue1}{rgb}{0.79,0.88,1.00}
\definecolor{lightsteelblue2}{rgb}{0.74,0.82,0.93}
\definecolor{lightsteelblue3}{rgb}{0.64,0.71,0.80}
\definecolor{lightsteelblue4}{rgb}{0.43,0.48,0.55}
\definecolor{lightyellow}{rgb}{1.00,1.00,0.88}
\definecolor{lightsteel}{rgb}{0.69,0.77,0.87}
\definecolor{lightyellow}{rgb}{1.00,1.00,0.88}
\definecolor{lightyellow1}{rgb}{1.00,1.00,0.88}
\definecolor{lightyellow2}{rgb}{0.93,0.93,0.82}
\definecolor{lightyellow3}{rgb}{0.80,0.80,0.71}
\definecolor{Lightyellow4}{rgb}{0.55,0.55,0.48}
\definecolor{limegreen}{rgb}{0.20,0.80,0.20}
\definecolor{linen}{rgb}{0.98,0.94,0.90}
\definecolor{magenta}{rgb}{1.00,0.00,1.00}
\definecolor{magenta1}{rgb}{1.00,0.00,1.00}
\definecolor{magenta2}{rgb}{0.93,0.00,0.93}
\definecolor{magenta3}{rgb}{0.80,0.00,0.80}
\definecolor{magenta4}{rgb}{0.55,0.00,0.55}
\definecolor{maroon}{rgb}{0.69,0.19,0.38}
\definecolor{maroon1}{rgb}{1.00,0.20,0.70}
\definecolor{maroon2}{rgb}{0.93,0.19,0.65}
\definecolor{maroon3}{rgb}{0.80,0.16,0.56}
\definecolor{maroon4}{rgb}{0.55,0.11,0.38}
\definecolor{mediumaquamarine}{rgb}{0.40,0.80,0.67}
\definecolor{mediumblue}{rgb}{0.00,0.00,0.80}
\definecolor{mediumorchid1}{rgb}{0.88,0.40,1.00}
\definecolor{mediumorchid2}{rgb}{0.82,0.37,0.93}
\definecolor{mediumorchid3}{rgb}{0.71,0.32,0.80}
\definecolor{mediumorchid4}{rgb}{0.48,0.22,0.55}
\definecolor{mediumorchid}{rgb}{0.73,0.33,0.83}
\definecolor{mediumpurple}{rgb}{0.58,0.44,0.86}
\definecolor{mediumpurple1}{rgb}{0.67,0.51,1.00}
\definecolor{mediumpurple2}{rgb}{0.62,0.47,0.93}
\definecolor{mediumpurple3}{rgb}{0.54,0.41,0.80}
\definecolor{mediumpurple4}{rgb}{0.36,0.28,0.55}
\definecolor{medium}{rgb}{0.45,0.45,0.45}
\definecolor{mediumseagreen}{rgb}{0.24,0.70,0.44}
\definecolor{mediumsea}{rgb}{0.24,0.70,0.44}
\definecolor{mediumslateblue}{rgb}{0.48,0.41,0.93}
\definecolor{mediumslate}{rgb}{0.48,0.41,0.93}
\definecolor{mediumspringgreen}{rgb}{0.00,0.98,0.60}
\definecolor{mediumspring}{rgb}{0.00,0.98,0.60}
\definecolor{mediumturquoise}{rgb}{0.28,0.82,0.80}
\definecolor{mediumvioletred}{rgb}{0.78,0.08,0.52}
\definecolor{mediumviolet}{rgb}{0.78,0.08,0.52}
\definecolor{midnightblue}{rgb}{0.10,0.10,0.44}
\definecolor{mintcream}{rgb}{0.96,1.00,0.98}
\definecolor{mistyrose}{rgb}{1.00,0.89,0.88}
\definecolor{mistyrose1}{rgb}{1.00,0.89,0.88}
\definecolor{mistyrose2}{rgb}{0.93,0.84,0.82}
\definecolor{mistyrose3}{rgb}{0.80,0.72,0.71}
\definecolor{mistyrose4}{rgb}{0.55,0.49,0.48}
\definecolor{moccasin}{rgb}{1.00,0.89,0.71}
\definecolor{navajowhite}{rgb}{1.00,0.87,0.68}
\definecolor{navajowhite1}{rgb}{1.00,0.87,0.68}
\definecolor{navajowhite2}{rgb}{0.93,0.81,0.63}
\definecolor{navajowhite3}{rgb}{0.80,0.70,0.55}
\definecolor{navajowhite4}{rgb}{0.55,0.47,0.37}
\definecolor{navyblue}{rgb}{0.00,0.00,0.50}
\definecolor{navy}{rgb}{0.00,0.00,0.50}
\definecolor{oldlace}{rgb}{0.99,0.96,0.90}
\definecolor{olivedrab}{rgb}{0.42,0.56,0.14}
\definecolor{olivedrab1}{rgb}{0.75,1.00,0.24}
\definecolor{olivedrab2}{rgb}{0.70,0.93,0.23}
\definecolor{olivedrab3}{rgb}{0.60,0.80,0.20}
\definecolor{olivedrab4}{rgb}{0.41,0.55,0.13}
\definecolor{orange1}{rgb}{1.00,0.65,0.00}
\definecolor{orange2}{rgb}{0.93,0.60,0.00}
\definecolor{orange3}{rgb}{0.80,0.52,0.00}
\definecolor{orange4}{rgb}{0.55,0.35,0.00}
\definecolor{orangered}{rgb}{1.00,0.27,0.00}
\definecolor{orangered1}{rgb}{1.00,0.27,0.00}
\definecolor{orangered2}{rgb}{0.93,0.25,0.00}
\definecolor{orangered3}{rgb}{0.80,0.22,0.00}
\definecolor{orangered4}{rgb}{0.55,0.15,0.00}
\definecolor{orange}{rgb}{1.00,0.65,0.00}
\definecolor{orchid}{rgb}{0.85,0.44,0.84}
\definecolor{orchid1}{rgb}{1.00,0.51,0.98}
\definecolor{orchid2}{rgb}{0.93,0.48,0.91}
\definecolor{orchid3}{rgb}{0.80,0.41,0.79}
\definecolor{orchid4}{rgb}{0.55,0.28,0.54}
\definecolor{palegoldenrod}{rgb}{0.93,0.91,0.67}
\definecolor{palegreen}{rgb}{0.60,0.98,0.60}
\definecolor{palegreen1}{rgb}{0.60,1.00,0.60}
\definecolor{palegreen2}{rgb}{0.56,0.93,0.56}
\definecolor{palegreen3}{rgb}{0.49,0.80,0.49}
\definecolor{palegreen4}{rgb}{0.33,0.55,0.33}
\definecolor{paleturquoise1}{rgb}{0.73,1.00,1.00}
\definecolor{paleturquoise2}{rgb}{0.68,0.93,0.93}
\definecolor{paleturquoise3}{rgb}{0.59,0.80,0.80}
\definecolor{paleturquoise4}{rgb}{0.40,0.55,0.55}
\definecolor{paleturquoise}{rgb}{0.69,0.93,0.93}
\definecolor{palevioletred}{rgb}{0.86,0.44,0.58}
\definecolor{palevioletred1}{rgb}{1.00,0.51,0.67}
\definecolor{palevioletred2}{rgb}{0.93,0.47,0.62}
\definecolor{palevioletred3}{rgb}{0.80,0.41,0.54}
\definecolor{palevioletred4}{rgb}{0.55,0.28,0.36}
\definecolor{paleviolet}{rgb}{0.86,0.44,0.58}
\definecolor{papayawhip}{rgb}{1.00,0.94,0.84}
\definecolor{peachPuff1}{rgb}{1.00,0.85,0.73}
\definecolor{peachPuff2}{rgb}{0.93,0.80,0.68}
\definecolor{peachPuff3}{rgb}{0.80,0.69,0.58}
\definecolor{peachPuff4}{rgb}{0.55,0.47,0.40}
\definecolor{peachpuff}{rgb}{1.00,0.85,0.73}
\definecolor{peru}{rgb}{0.80,0.52,0.25}
\definecolor{pink}{rgb}{1.00,0.75,0.80}
\definecolor{pink1}{rgb}{1.00,0.71,0.77}
\definecolor{pink2}{rgb}{0.93,0.66,0.72}
\definecolor{pink3}{rgb}{0.80,0.57,0.62}
\definecolor{pink4}{rgb}{0.55,0.39,0.42}
\definecolor{plum}{rgb}{0.87,0.63,0.87}
\definecolor{plum1}{rgb}{1.00,0.73,1.00}
\definecolor{plum2}{rgb}{0.93,0.68,0.93}
\definecolor{plum3}{rgb}{0.80,0.59,0.80}
\definecolor{plum4}{rgb}{0.55,0.40,0.55}
\definecolor{powderblue}{rgb}{0.69,0.88,0.90}
\definecolor{purple}{rgb}{0.63,0.13,0.94}
\definecolor{purple1}{rgb}{0.61,0.19,1.00}
\definecolor{purple2}{rgb}{0.57,0.17,0.93}
\definecolor{purple3}{rgb}{0.49,0.15,0.80}
\definecolor{purple4}{rgb}{0.33,0.10,0.55}
\definecolor{red}{rgb}{1.00,0.00,0.00}
\definecolor{red1}{rgb}{1.00,0.00,0.00}
\definecolor{red2}{rgb}{0.93,0.00,0.00}
\definecolor{red3}{rgb}{0.80,0.00,0.00}
\definecolor{red4}{rgb}{0.55,0.00,0.00}
\definecolor{rosybrown}{rgb}{0.74,0.56,0.56}
\definecolor{rosybrown1}{rgb}{1.00,0.76,0.76}
\definecolor{rosybrown2}{rgb}{0.93,0.71,0.71}
\definecolor{rosybrown3}{rgb}{0.80,0.61,0.61}
\definecolor{rosybrown4}{rgb}{0.55,0.41,0.41}
\definecolor{royalblue}{rgb}{0.25,0.41,0.88}
\definecolor{royalblue1}{rgb}{0.28,0.46,1.00}
\definecolor{royalblue2}{rgb}{0.26,0.43,0.93}
\definecolor{royalblue3}{rgb}{0.23,0.37,0.80}
\definecolor{royalblue4}{rgb}{0.15,0.25,0.55}
\definecolor{saddlebrown}{rgb}{0.55,0.27,0.07}
\definecolor{salmon}{rgb}{0.98,0.50,0.45}
\definecolor{salmon1}{rgb}{1.00,0.55,0.41}
\definecolor{salmon2}{rgb}{0.93,0.51,0.38}
\definecolor{salmon3}{rgb}{0.80,0.44,0.33}
\definecolor{salmon4}{rgb}{0.55,0.30,0.22}
\definecolor{sandybrown}{rgb}{0.96,0.64,0.38}
\definecolor{seagreen}{rgb}{0.18,0.55,0.34}
\definecolor{seagreen1}{rgb}{0.33,1.00,0.62}
\definecolor{seagreen2}{rgb}{0.31,0.93,0.58}
\definecolor{seagreen3}{rgb}{0.26,0.80,0.50}
\definecolor{seagreen4}{rgb}{0.18,0.55,0.34}
\definecolor{seashell}{rgb}{1.00,0.96,0.93}
\definecolor{seashell1}{rgb}{1.00,0.96,0.93}
\definecolor{seashell2}{rgb}{0.93,0.90,0.87}
\definecolor{seashell3}{rgb}{0.80,0.77,0.75}
\definecolor{seashell4}{rgb}{0.55,0.53,0.51}
\definecolor{sienna}{rgb}{0.63,0.32,0.18}
\definecolor{sienna1}{rgb}{1.00,0.51,0.28}
\definecolor{sienna2}{rgb}{0.93,0.47,0.26}
\definecolor{sienna3}{rgb}{0.80,0.41,0.22}
\definecolor{sienna4}{rgb}{0.55,0.28,0.15}
\definecolor{skyblue1}{rgb}{0.53,0.81,1.00}
\definecolor{skyblue2}{rgb}{0.49,0.75,0.93}
\definecolor{skyblue3}{rgb}{0.42,0.65,0.80}
\definecolor{skyblue4}{rgb}{0.29,0.44,0.55}
\definecolor{skyblue}{rgb}{0.53,0.81,0.92}
\definecolor{slateblue1}{rgb}{0.51,0.44,1.00}
\definecolor{slateblue2}{rgb}{0.48,0.40,0.93}
\definecolor{slateblue3}{rgb}{0.41,0.35,0.80}
\definecolor{slateblue4}{rgb}{0.28,0.24,0.55}
\definecolor{slateblue}{rgb}{0.42,0.35,0.80}
\definecolor{slategray}{rgb}{0.44,0.50,0.56}
\definecolor{slategray1}{rgb}{0.78,0.89,1.00}
\definecolor{slategray2}{rgb}{0.73,0.83,0.93}
\definecolor{slategray3}{rgb}{0.62,0.71,0.80}
\definecolor{slategray4}{rgb}{0.42,0.48,0.55}
\definecolor{snow}{rgb}{1.00,0.98,0.98}
\definecolor{snow1}{rgb}{1.00,0.98,0.98}
\definecolor{snow2}{rgb}{0.93,0.91,0.91}
\definecolor{snow3}{rgb}{0.80,0.79,0.79}
\definecolor{snow4}{rgb}{0.55,0.54,0.54}
\definecolor{springgreen1}{rgb}{0.00,1.00,0.50}
\definecolor{springgreen2}{rgb}{0.00,0.93,0.46}
\definecolor{springgreen3}{rgb}{0.00,0.80,0.40}
\definecolor{springgreen4}{rgb}{0.00,0.55,0.27}
\definecolor{springgreen}{rgb}{0.00,1.00,0.50}
\definecolor{steelblue}{rgb}{0.27,0.51,0.71}
\definecolor{steelblue1}{rgb}{0.39,0.72,1.00}
\definecolor{steelblue2}{rgb}{0.36,0.67,0.93}
\definecolor{steelblue3}{rgb}{0.31,0.58,0.80}
\definecolor{steelblue4}{rgb}{0.21,0.39,0.55}
\definecolor{tan}{rgb}{0.82,0.71,0.55}
\definecolor{tan1}{rgb}{1.00,0.65,0.31}
\definecolor{tan2}{rgb}{0.93,0.60,0.29}
\definecolor{tan3}{rgb}{0.80,0.52,0.25}
\definecolor{tan4}{rgb}{0.55,0.35,0.17}
\definecolor{thistle}{rgb}{0.85,0.75,0.85}
\definecolor{thistle1}{rgb}{1.00,0.88,1.00}
\definecolor{thistle2}{rgb}{0.93,0.82,0.93}
\definecolor{thistle3}{rgb}{0.80,0.71,0.80}
\definecolor{thistle4}{rgb}{0.55,0.48,0.55}
\definecolor{tomato}{rgb}{1.00,0.39,0.28}
\definecolor{tomato1}{rgb}{1.00,0.39,0.28}
\definecolor{tomato2}{rgb}{0.93,0.36,0.26}
\definecolor{tomato3}{rgb}{0.80,0.31,0.22}
\definecolor{tomato4}{rgb}{0.55,0.21,0.15}
\definecolor{turquoise}{rgb}{0.25,0.88,0.82}
\definecolor{turquoise1}{rgb}{0.00,0.96,1.00}
\definecolor{turquoise2}{rgb}{0.00,0.90,0.93}
\definecolor{turquoise3}{rgb}{0.00,0.77,0.80}
\definecolor{turquoise4}{rgb}{0.00,0.53,0.55}
\definecolor{violetred1}{rgb}{1.00,0.24,0.59}
\definecolor{violetred2}{rgb}{0.93,0.23,0.55}
\definecolor{violetred3}{rgb}{0.80,0.20,0.47}
\definecolor{violetred4}{rgb}{0.55,0.13,0.32}
\definecolor{violetred}{rgb}{0.82,0.13,0.56}
\definecolor{violet}{rgb}{0.93,0.51,0.93}
\definecolor{wheat}{rgb}{0.96,0.87,0.70}
\definecolor{wheat1}{rgb}{1.00,0.91,0.73}
\definecolor{wheat2}{rgb}{0.93,0.85,0.68}
\definecolor{wheat3}{rgb}{0.80,0.73,0.59}
\definecolor{wheat4}{rgb}{0.55,0.49,0.40}
\definecolor{white}{rgb}{1.00,1.00,1.00}
\definecolor{whitesmoke}{rgb}{0.96,0.96,0.96}
\definecolor{yellow}{rgb}{1.00,1.00,0.00}
\definecolor{yellow1}{rgb}{1.00,1.00,0.00}
\definecolor{yellow2}{rgb}{0.93,0.93,0.00}
\definecolor{yellow3}{rgb}{0.80,0.80,0.00}
\definecolor{yellow4}{rgb}{0.55,0.55,0.00}
\definecolor{yellowgreen}{rgb}{0.60,0.80,0.20}
\makeatletter
\newlength\@tempdim@x
\newlength\@tempdim@y
\newcommand\AtUpperLeftCorner[3]{%
\begingroup
\@tempdim@x=0cm
\@tempdim@y=\paperheight
\advance\@tempdim@x#1
\advance\@tempdim@y-#2
\put(\LenToUnit{\@tempdim@x},\LenToUnit{\@tempdim@y}){#3}%
\endgroup}
\newcommand\AtUpperRightCorner[3]{%
\begingroup
\@tempdim@x=\paperwidth
\@tempdim@y=\paperheight
\advance\@tempdim@x-#1
\advance\@tempdim@y-#2
\put(\LenToUnit{\@tempdim@x},\LenToUnit{\@tempdim@y}){#3}%
\endgroup}
\newcommand\AtLowerLeftCorner[3]{%
\begingroup
\@tempdim@x=0cm
\@tempdim@y=0cm
\advance\@tempdim@x#1
\advance\@tempdim@y#2
\put(\LenToUnit{\@tempdim@x},\LenToUnit{\@tempdim@y}){#3}%
\endgroup}
\newcommand\AtLowerRightCorner[3]{%
\begingroup
\@tempdim@x=\paperwidth
\@tempdim@y=0cm
\advance\@tempdim@x-#1
\advance\@tempdim@y#2
\put(\LenToUnit{\@tempdim@x},\LenToUnit{\@tempdim@y}){#3}%
\endgroup}
\makeatother
\newtheorem{theoreme}{Theorem}[section]

\newtheorem{definition}[theoreme]{Definition}
\newtheorem{proposition}[theoreme]{Proposition}
\newtheorem{lemme}[theoreme]{Lemma}
\newtheorem{corollaire}[theoreme]{Corollary}

\newtheorem{remarque}[theoreme]{Remark}
\newtheorem{remarques}[theoreme]{Remarks}
\newenvironment{demo}{\begin{proof}}{\end{proof}}

\addto\captionsfrenchb{}
\addto\captionsfrenchb{}

    \newlength{\myarrowsize} 
    \newlength{\myoldlinewidth}

    \pgfarrowsdeclare{myto}{myto}{
        \pgfsetdash{}{0pt} 
        \pgfsetbeveljoin 
        \pgfsetroundcap 
        \setlength{\myarrowsize}{0.5pt}
        \addtolength{\myarrowsize}{.5\pgflinewidth}
        \pgfarrowsleftextend{-4\myarrowsize-.5\pgflinewidth} 
        \pgfarrowsrightextend{.7\pgflinewidth}
    }{
        \setlength{\myarrowsize}{0.4pt} 
        \addtolength{\myarrowsize}{.3\pgflinewidth}  
        \setlength{\myoldlinewidth}{\pgflinewidth}
        \pgfsetroundjoin
        \pgfsetlinewidth{0.0001pt}
        \pgfpathmoveto{\pgfpoint{0.43\myarrowsize}{0}}
        \pgfpatharc{0}{70}{0.14\myarrowsize}
        \pgfpatharc{-80}{-169.5}{4\myarrowsize}
        \pgfpatharc{150}{189}{0.95\myarrowsize and 0.95\myarrowsize}
        \pgfpatharc{0}{-40}{15\myarrowsize}
        \pgfpathmoveto{\pgfpoint{0.43\myarrowsize}{0}}
        \pgfpatharc{0}{-70}{0.14\myarrowsize}
        \pgfpatharc{80}{169.5}{4\myarrowsize}
        \pgfpatharc{-150}{-189}{0.95\myarrowsize and 0.95\myarrowsize}
        \pgfpatharc{0}{40}{15\myarrowsize}
        \pgfpathclose
        \pgfsetstrokeopacity{0.25}
        \pgfusepathqfillstroke
    }

    \pgfarrowsdeclare{myonto}{myonto}{
        \pgfsetdash{}{0pt} 
        \pgfsetbeveljoin 
        \pgfsetroundcap 
        \setlength{\myarrowsize}{0.5pt}
        \addtolength{\myarrowsize}{.5\pgflinewidth}
        \pgfarrowsleftextend{-4\myarrowsize-.5\pgflinewidth} 
        \pgfarrowsrightextend{.7\pgflinewidth}
    }{
        \setlength{\myarrowsize}{0.4pt} 
        \addtolength{\myarrowsize}{.3\pgflinewidth}  
        \setlength{\myoldlinewidth}{\pgflinewidth}
        \pgfsetroundjoin
        \pgfsetlinewidth{0.0001pt}
        \pgfpathmoveto{\pgfpoint{0.43\myarrowsize}{0}}
        \pgfpatharc{0}{70}{0.14\myarrowsize}
        \pgfpatharc{-80}{-169.5}{4\myarrowsize}
        \pgfpatharc{150}{189}{0.95\myarrowsize and 0.95\myarrowsize}
        \pgfpatharc{0}{-40}{15\myarrowsize}
        \pgfpathmoveto{\pgfpoint{-7\myarrowsize}{0}}
				\pgfpatharc{0}{70}{0.14\myarrowsize}
        \pgfpatharc{-80}{-169.5}{4\myarrowsize}
        \pgfpatharc{150}{189}{0.95\myarrowsize and 0.95\myarrowsize}
        \pgfpatharc{0}{-40}{15\myarrowsize}
        \pgfpathmoveto{\pgfpoint{0.43\myarrowsize}{0}}
        \pgfpatharc{0}{-70}{0.14\myarrowsize}
        \pgfpatharc{80}{169.5}{4\myarrowsize}
        \pgfpatharc{-150}{-189}{0.95\myarrowsize and 0.95\myarrowsize}
        \pgfpatharc{0}{40}{15\myarrowsize}
			  \pgfpathmoveto{\pgfpoint{-7\myarrowsize}{0}}
        \pgfpatharc{0}{-70}{0.14\myarrowsize}
        \pgfpatharc{80}{169.5}{4\myarrowsize}
        \pgfpatharc{-150}{-189}{0.95\myarrowsize and 0.95\myarrowsize}
        \pgfpatharc{0}{40}{15\myarrowsize}
        \pgfpathclose
        \pgfsetstrokeopacity{0.25}
        \pgfusepathqfillstroke
        \pgfusepathqfillstroke
    }

    \pgfarrowsdeclare{myhook}{myhook}{
        \setlength{\myarrowsize}{0.6pt}
        \addtolength{\myarrowsize}{.5\pgflinewidth}
        \pgfarrowsleftextend{-4\myarrowsize-.5\pgflinewidth} 
        \pgfarrowsrightextend{.7\pgflinewidth}
    }{
        \setlength{\myarrowsize}{0.6pt} 
        \addtolength{\myarrowsize}{.5\pgflinewidth}  
        \pgfsetdash{}{+0pt}
        \pgfsetroundcap
        \pgfpathmoveto{\pgfqpoint{-2pt}{-6\pgflinewidth}}
        \pgfpathcurveto
            {\pgfqpoint{4\pgflinewidth}{-4.667\pgflinewidth}}
            {\pgfqpoint{4\pgflinewidth}{0pt}}
            {\pgfpointorigin}
        \pgfusepathqstroke
    }

		\pgfarrowsdeclare{my to}{my to}
{
  \pgfarrowsleftextend{-2\pgflinewidth}
  \pgfarrowsrightextend{\pgflinewidth}
}
{
  \pgfsetlinewidth{0.8\pgflinewidth}
  \pgfsetdash{}{0pt}
  \pgfsetroundcap
  \pgfsetroundjoin
  \pgfpathmoveto{\pgfpoint{-5.5\pgflinewidth}{7.5\pgflinewidth}}
  \pgfpathcurveto
  {\pgfpoint{-4.0\pgflinewidth}{0.1\pgflinewidth}}
  {\pgfpoint{0pt}{0.25\pgflinewidth}}
  {\pgfpoint{0.75\pgflinewidth}{0pt}}
  \pgfpathcurveto
  {\pgfpoint{0pt}{-0.25\pgflinewidth}}
  {\pgfpoint{-4.0\pgflinewidth}{-0.1\pgflinewidth}}
  {\pgfpoint{-5.5\pgflinewidth}{-7.5\pgflinewidth}}
  \pgfusepathqstroke
}

\tikzstyle{vecArrow} = [thick, decoration={markings,mark=at position
   1 with {\arrow[semithick]{open triangle 60}}},
   double distance=1.4pt, shorten >= 5.5pt,
   preaction = {decorate},
   postaction = {draw,line width=1.4pt, white,shorten >= 4.5pt}]
\tikzstyle{innerWhite} = [semithick, white,line width=1.4pt, shorten >= 4.5pt]

	\makeatletter
	\newcommand\POSITION[3]{%
	\begingroup
	\@tempdim@x=0cm
	\@tempdim@y=\paperheight
	\advance\@tempdim@x#1
	\advance\@tempdim@y-#2
	\put(\LenToUnit{\@tempdim@x},\LenToUnit{\@tempdim@y}){#3}%
	\endgroup
	}
\makeatother

\addto\captionsenglish{}

\begin{document}
	\maketitle
	
	\begin{abstract}
		We classify fully commutative elements in the affine Coxeter group of type $\tilde{A_{n}}$. We give a normal form for such elements, then we propose an application of this normal form: we lift these fully commutative elements to the affine braid group of type $\tilde{A_{n}}$ and we get a form for ``fully commutative braids''.   
	\end{abstract}
	
	
	\section{Introduction and notation}
		
		 	Let $(W,S)$ be a Coxeter system with associated Dynkin Diagram $\Gamma$. Let $w\in W$. We know that from a given reduced expression of $w$ we can arrive to any other reduced expression only by applying braid relations \cite{Bourbaki_1981}. Among these relations there are commutation relations corresponding to the non-neighbours (precisely, $t$ and $s$ with $m_{st} = 2$). \\
			
			\begin{definition}
			Let $w$ be in $W$. We call support of $w$ the subset of $S$ consisting of all generators appearing in a (any) reduced expression of $w$ . It is to be denoted by $Supp(w)$.  \\
		\end{definition}
		
		We define $\mathscr{L} (w) $ to be the set of $s\in S$ such that $l(sw)<l(w)$, in other terms  $s$ appears at the left edge of some reduced expression of $w$. Similarly we define $\mathscr{R}(w)$.\\

			\begin{definition}
			Elements for which one can pass from any reduced expression to any other one only by applying commutation relations are called fully commutative elements. Usually we denote the set of fully commutative elements by $W^{c}$. \\
		    \end{definition}
			
			Consider the $A$-type Coxeter group with $n$ generators $W(A_{n})$, with the following Dynkin diagram:\\
			
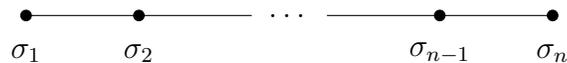
			\begin{figure}[ht]
				\centering
				\begin{tikzpicture}

  \filldraw (0,0) circle (2pt);
  \node at (0,-0.5) {$\sigma_{1}$}; 
   
  \draw (0,0) -- (1.5, 0);

  \filldraw (1.5,0) circle (2pt);
  \node at (1.5,-0.5) {$\sigma_{2}$};

  \draw (1.5,0) -- (3, 0);

  \node at (3.5,0) {$\dots$};

  \draw (4,0) -- (5.5, 0);
  
  \filldraw (5.5,0) circle (2pt);
  \node at (5.5,-0.5) {$\sigma_{n-1}$};
 
  \draw (5.5,0) -- (7, 0);
  
  \filldraw (7,0) circle (2pt);
  \node at (7,-0.5) {$\sigma_{n}$};

               \end{tikzpicture}
			\end{figure}
			
			We set $W^{c}(A_{n})$ to be the set of its fully commutative elements, its cardinality is the Catalan number $\frac{1}{n+2}\big(^{~2(n+1)~}_{~~n+1}\big)$. One can prove easily by induction on $n$ $\big($considering right classes of $W(A_{n-1})$ in $W(A_{n})\big)$ the following well known theorem.\\
			
			\begin{theoreme} \label{1_2}
				Let $u$ be any fully commutative element in $W(A_{n})$. Then there is a unique reduced expression of $u$ of the form:
				\begin{eqnarray}
					u = \sigma_{i_{1}} \sigma_{i_{1}-1} ... \sigma_{j_{1}}~~ \sigma_{i_{2}} \sigma_{i_{2}-1}... \sigma_{j_{2}} ~~...~~\sigma_{i_{p}}\sigma_{i_{p}-1} ... \sigma_{j_{p}}, \nonumber
				\end{eqnarray}
				
				where $ 1\leq i_{1} < i_{2} .. < i_{p} \leq n $, $ 1\leq j_{1} < j_{2} .. < j_{p} \leq n $ and $ j_{k} \leq i_{k}$ for every $1 \leq k \leq p $. \\ 
			\end{theoreme}			
			
			Now let $W(\tilde{A_{n}}) $ be the affine Coxeter group of $\tilde{A}$-type with $n+1$ generators. with the following Dynkin diagram: 
			
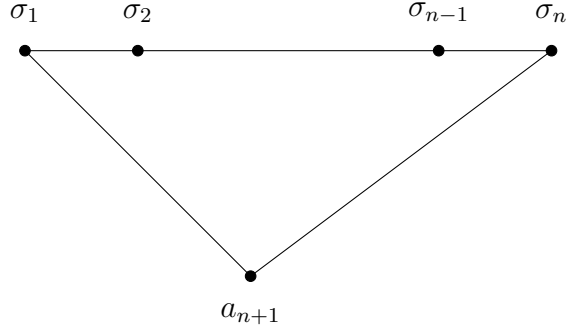
\begin{figure}[ht]
				\centering
				\begin{tikzpicture}

 \node at (0,0.5) {$\sigma_{1}$}; 
  \filldraw (0,0) circle (2pt);
   
  \draw (0,0) -- (1.5, 0);
  
  \node at (1.5,0.5) {$\sigma_{2}$};
  \filldraw (1.5,0) circle (2pt);

  \draw (1.5,0) -- (5.5, 0);

  \node at (5.5,0.5) {$\sigma_{n-1}$};
  \filldraw (5.5,0) circle (2pt);
 
  \draw (5.5,0) -- (7, 0);
  
  \node at (7,0.5) {$\sigma_{n}$};
  \filldraw (7,0) circle (2pt);

  \draw (7,0) -- (3, -3);
  
  \filldraw (3, -3) circle (2pt);
  \node at (3, -3.5) {$a_{n+1}$};

  \draw (3, -3) -- (0, 0);
               \end{tikzpicture}
			\end{figure}
			
			The main result of this paper is presented in theorem \ref{2_4_4} which is the affine version of theorem \ref{1_2}.\\
			
			In the second section, we give some general definitions. Then we state and prove our general result about the affine fully commutative elements. In the third section we give a consequence of our classification. We lift the fully commutative elements to elements having the same expressions in the $A$-type braid group  $B(\tilde{A_{n}}) $, or : ``fully commutative braids'' which will be the key, in a forthcoming paper, to classify traces on the associated affine Temperley-Lieb algebra\cite{Sadek_2013_2}.\\

	\section{Fully commutative elements} 
		
			Let $(W,S)$ be a Coxeter system such that any two elements in $S$ are conjugate in $W$, in this case fully commutative elements have some additional elegant properties, for example we can reformulate the definition as follows.\\
			
		\begin{proposition}
			Let $(W,S)$ be such that any two elements in $S$ are conjugate in $W$. Let $w\in W$. Then $w$ is fully commutative if and only if every $s$ in $Supp(w)$ occurs the same number of times in any reduced expression of $w$. \\
		\end{proposition}
		
		\begin{demo}
		We omit the proof. 
		\end{demo}
			
		Hence, in this case, for a fully commutative element $w$, we can talk of the multiplicity of a simple reflexion in $Supp(w)$. That is if $s$ is in $Supp(w)$, we call the multiplicity of $s$ in $w$ the number of times $s$ appears in a (hence every) reduced expression of $w$. The center of our interest in this work is fully commutative elements in $\tilde{A}$-type Coxeter groups, which is an example of Coxeter groups in which any two elements in $S$ are conjugate.\\
		
	    Notice that, in theorem \ref{1_2}, if $ \sigma_{n} $ belongs to $ supp (u)$, then  $ \sigma_{n} $  will certainly appear only once, and it is to be equal to $\sigma_{i_{p}}$. Similarly for $\sigma_{1} $: if it belongs to $ supp (u)$, then $\sigma_{1}$ will certainly appear only once, and it is equal to $\sigma_{j_{1}}$.\\
			
			\begin{definition}
				An element $u$ in $W^{c}(A_{n})$ is called full if and only if both $ \sigma_{n} $ and $ \sigma_{1} $ belong to $Supp(w)$. In this case $u$ has a reduced expression of the form:  
				\begin{eqnarray}
					u = \sigma_{i_{1}} .. \sigma_{1} \sigma_{i_{2}} .. \sigma_{j_{2}} ~...~ \sigma_{n} .. \sigma_{j_{p}},\nonumber
				\end{eqnarray}
				
				where $ 1\leq i_{1} < i_{2} .. < i_{p-1} \leq n $, $ 1\leq j_{2} .. < j_{p} \leq n $ and $ j_{k} \leq i_{k}$ for every $1\leq k \leq n $.\\
			\end{definition}
			
			\begin{definition}
				Suppose that  $u$ is full, i.e., $ u = \sigma_{i_{1}} .. \sigma_{1} \sigma_{i_{2}} .. \sigma_{j_{2}} ~...~ \sigma_{n} .. \sigma_{j_{p}} $. We say that $\sigma_{n}$ is on the left (in $u$), if and only if $ u = \sigma_{n} ~...~  \sigma_{2}  \sigma_{1} $. In all other cases we say that $\sigma_{n} $ is on the right.\\
			\end{definition}

		\subsection{Classification of $W^{c}(\tilde{A_{n}}) $: a normal form}
			\vspace{0.25cm}

	        In this subsection we prove the following theorem.
			
			\begin{theoreme} \label{2_4_4}
				Let $ 2 \leq n $. Let $w \in W(\tilde{A_{n}})$ be a fully commutative such that $a_{n+1} \in supp(w)$. Then, there exists a unique reduced expression of $w$, of the following form:\\
				\vspace{-0.25cm}
				\begin{eqnarray}
					w &=& \sigma_{i_{1}} .. \sigma_{2} \sigma_{1} \sigma_{r_{1}} .. \sigma_{n-1} \sigma_{n}  a_{n+1} \sigma_{i_{2}} .. \sigma_{2} \sigma_{1} \sigma_{r_{2}} .. \sigma_{n-1} \sigma_{n}  a_{n+1} ~...~ \sigma_{i_{p}} .. \sigma_{2} \sigma_{1} \sigma_{r_{p}} .. \sigma_{n-1} \sigma_{n}\nonumber\\
					& & (a_{n+1}  \sigma_{j} ..\sigma_{2} \sigma_{1}  \sigma_{j +1 } ..\sigma_{n-1} \sigma_{n})^{k} u ,\nonumber
				\end{eqnarray}
				
				where: $ 0 \leq i_{1} < i_{2}~...~ <i_{p} < r_{p} < r_{2} ~...~< r_{1} \leq n +1  $, $ r_{p} - i_{p} \geq 2 $, $ i_{p} <  j  \leq r_{p} -1 $, $i_{1} \leq n $ and $ 0 \leq k $, \\ 
	
				and where $u$ has one of the following forms\\ 
				
				\begin{itemize}[label=$\bullet$, font=\normalsize, font=\color{black}, leftmargin=2cm,parsep=0cm, itemsep=0.25cm, topsep=0cm]
					\item If $k=0$, then:
						\vspace{-0.25cm}
						\begin{eqnarray}
							u= a_{n+1} \sigma_{l_{1}} .. \sigma_{g_{1}} \sigma_{l_{2}} .. \sigma_{g_{2}} ~...~ \sigma_{l_{t}} .. \sigma_{g_{t}}, \nonumber
						\end{eqnarray}
						\vspace{0.25cm}
						where $ 1\leq l_{1} < l_{2} .. < l_{t} \leq n $, $ 1\leq g_{1} < g_{2} .. < g_{t} \leq n $ and $ g_{k} \leq l_{k}$, for any $1 \leq k \leq n $. With $ i_{p} < l_{1} $ and $g_{t} <r_{p} $.\\ 
						
					\item If $k \geq 1 $, then:
						\vspace{-0.25cm}
						\begin{eqnarray}
							u = a_{n+1} \sigma_{j} ~...~ \sigma_{d_{1}} \sigma_{j +1 } ~...~ \sigma_{d_{2}} \sigma_{j +2} ~...~ \sigma_{d_{3}} ~...~ \sigma_{j + z} ~...~ \sigma_{d_{z+1}}, \nonumber
						\end{eqnarray}
				
						where $ d_{1} < d_{2} ~...~  < d_{z+1} $ and $ j + c \geq d_{ c+1} $ for  $ 0 \leq c \leq z $.  \\
				\end{itemize}
			\end{theoreme}
			\begin{definition}   
				We define the affine length of $u$ in $W^{c}(A_{n})$ to be the multiplicity of $a_{n+1}$ in $Supp(u)$. We denote it by $L(u)$. \\
			\end{definition}
			
			Suppose that $w$ is a fully commutative element in $W(\tilde{A_{n}})$. Clearly $ L (w) = 0 $ expresses the case where $a_{n+1}$ is not in $ supp (w)$, in other terms $ w $ is a fully commutative element in  $W(A_{n})$. Suppose that $L(w)= m$ where $m$ is positive. Any reduced expression of $w$ is of the form: 
			\vspace{-0.25cm}
			\begin{eqnarray}
				w = u_{1} a_{n+1} u_{2} a_{n+1} ~...~ u_{m} a_{n+1} u_{m+1}, \nonumber
			\end{eqnarray}
			
			where $ u_{i} $ is in $W^{c}(A_{n})$, for $1 \leq i \leq m+1 $. Moreover, suppose that $ L (w) \geq 2 $. Then $u_{i}$ must be full for $ 2\leq u_{i} \leq m $, otherwise $w$ is not fully commutative. \\
	
	    		Before treating the general case, we classify fully commutative elements of $W(\tilde{A_{2}})$. This gives an idea about the general proof, in its simplest form.\\
			
			\begin{theoreme}\label{2_6}
				Let $w$ be in $W^{c}(\tilde{A_{2}})$. Then there exists $0 \leq k$, such that  $w$ has one and only one of the following forms:\\ 
				
				\vspace{0.5cm}
				\begin{tikzpicture}
               \begin{scope}[xscale = 1]

  \node[left =-13pt] at (-1.5,1)  {$1$};
  \node[left =-13pt] at (-1.5,0)  {$a_{3}$};
  \node[left =-13pt] at (-1.5,-1) {$\sigma_{1} a_{3}$};

	\draw (-1,1)  -- (0,0);
	\draw (-1,0)  -- (0,0);
	\draw (-1,-1) -- (0,0);
	
  \node at (1,0)  {$(\sigma_{2}\sigma_{1}a_{3})^{k}$};

	\draw (2,0)  -- (3,-1);
	\draw (2,0)  -- (3,0);
	\draw (2,0)  -- (3,1);

  \node[right =-13pt] at (3.5,1)  {$1$};
  \node[right =-13pt] at (3.5,0)  {$\sigma_{2}$};
  \node[right =-13pt] at (3.5,-1) {$\sigma_{2} \sigma_{1}$};
\end{scope}
\begin{scope}[xscale = 1, xshift = 8cm]

  \node[left =-13pt] at (-1.5,1)  {$1$};
  \node[left =-13pt] at (-1.5,0)  {$a_{3}$};
  \node[left =-13pt] at (-1.5,-1) {$\sigma_{2} a_{3}$};

	\draw (-1,1)  -- (0,0);
	\draw (-1,0)  -- (0,0);
	\draw (-1,-1) -- (0,0);
	
  \node at (1,0)  {$(\sigma_{1}\sigma_{2}a_{3})^{k}$};

	\draw (2,0)  -- (3,-1);
	\draw (2,0)  -- (3,0);
	\draw (2,0)  -- (3,1);

  \node[right =-13pt] at (3.5,1)  {$1$};
  \node[right =-13pt] at (3.5,0)  {$\sigma_{1}$};
  \node[right =-13pt] at (3.5,-1) {$\sigma_{1} \sigma_{2}$};

\end{scope}
               \end{tikzpicture}

			\end{theoreme}
			
			\vspace{0.4cm}
			
			\begin{demo}
			
				As we saw above $ w = u_{1} a_{3} u_{2} a_{3} ~...~ u_{m} a_{3} u_{m+1} $, where $u_{i}$ is in  $W^{c}(A_{2})$. If $L(w)$ is 0 or 1 it is clear that we can get it from the tree formulas above. Suppose that $2 \leq L(w)$. Hence $u_{i}$ is full for $ 2 \leq i \leq m $. In particular $u_{2}$ is full. Actually there are not many choices for $u_{2}$, since the only full elements in $W^{c}(A_{2})$ are $ \sigma_{1} \sigma_{2}$ and $ \sigma_{2} \sigma_{1}$. The first possibility is that $u_{2}=\sigma_{1} \sigma_{2}$. Now being a full element, $u_{3}$ is definitely equal to $\sigma_{1} \sigma_{2}$, otherwise we would have, in $w$, the following subword $ u_{1} a_{3}\sigma_{1} \underbrace{\sigma_{2} a_{3} \sigma_{2} }_{}\sigma_{1}$. This is not possible since $w$ is fully commutative, thus $u_{3} = u_{2} =  \sigma_{1} \sigma_{2} $. The same holds for every $u_{i}$ for $i \leq m$, i.e., if $u_{2}$ is equal to $ \sigma_{1} \sigma_{2}$ then $ w =  u_{1}a_{3}(\sigma_{1} \sigma_{2}  a_{3})^{m-1}u_{m+1} $.\\ 
			
				It is clear that $u_{1}$ is in $W^{c}(A_{2})$, and does not end with $\sigma_{1}$, hence $u_{1}$ is equal to $\sigma_{2}$ or 1. In the same way, we see that $u_{m+1}$ is in $W^{c}(A_{2})$, it cannot ends with $\sigma_{2}$, so $u_{m+1}$ is equal to $\sigma_{2}$,$\sigma_{1}\sigma_{2}$ or 1. In other terms, if $u_{2}$ is equal to $\sigma_{1}\sigma_{2}$ we get the second tree.\\
			
				Now suppose that $u_{2}=\sigma_{2}\sigma_{1}$, then $w =  u_{1}a_{3}(\sigma_{2}\sigma_{1} a_{3})^{m-1}u_{m+1}$. With a similar discussion about the first choice of $ u_{2}$, we see that when $u_{2}=\sigma_{2}\sigma_{1}$ we get the first tree. 
			
			\end{demo}
		
			In order to simplify, we suppose now that $ n \geq 3 $ (although many propositions in what follows are valid in $W(\tilde{A_{2}})$). \\
	    
			\begin{remarque}
				Let $u$ be a full element : $ u = \sigma_{i_{1}} .. \sigma_{1} \sigma_{i_{2}} .. \sigma_{j_{2}} ~...~ \sigma_{n} .. \sigma_{j_{p}} $. Assume that $ \sigma_{n} $ is on the right  in $u$, hence, by pushing $\sigma_{n} $ to the left  we see easily that 
				\begin{eqnarray}
					u = \sigma_{i} ..\sigma_{2} \sigma_{1}  \sigma_{r} ~...~\sigma_{n-1} \sigma_{n} x, \nonumber
				\end{eqnarray}
				
				where $ 1 \leq i \leq n-1 $, $ 1 \leq r \leq n $ and $  i < r $, while $ supp (x) \subseteq  \left\{ \sigma_{2},  \sigma_{3}.. \sigma_{n-1}   \right\} $ if $x$ is not $1$. \\ 
	
			\end{remarque}
			
			\begin{lemme}\label{2_4_11}
				Let $w$ be in $W^{c}(\tilde{A_{n}})$ such that $ 2 \leq L(w) $. Say: 
				\begin{eqnarray}
					w = u_{1} a_{n+1} u_{2} a_{n+1} ~...~ u_{m} a_{n+1} u_{m+1}. \nonumber
				\end{eqnarray}
				
				Assume that $\sigma_{n}$ is on the right in $ u_{h} $, for $ 2 \leq h \leq m$. Then $w$ has one of the three following forms:
				\begin{eqnarray}
					w_{1} &=& u_{1} a_{n+1} \sigma_{i_{1}} .. \sigma_{2} \sigma_{1} \sigma_{r_{1}} .. \sigma_{n-1} \sigma_{n} a_{n+1} \sigma_{i_{2}} .. \sigma_{2} \sigma_{1} \sigma_{r_{2}} .. \sigma_{n-1} \sigma_{n}  a_{n+1} ~...~ \sigma_{i_{p}} .. \sigma_{2} \sigma_{1} \sigma_{r_{p}} .. \sigma_{n-1} \sigma_{n} \nonumber\\
					& & (a_{n+1} \sigma_{j} .. \sigma_{2} \sigma_{1} \sigma_{j +1 } .. \sigma_{n-1} \sigma_{n})^{m-(1+p)} \nonumber\\
					& & a_{n+1} \sigma_{j} .. \sigma_{d_{1}} \sigma_{j +1} .. \sigma_{d_{2}} \sigma_{j +2} .. \sigma_{d_{3}} ~...~ \sigma_{j + z} ..\sigma_{d_{z+1}}, \nonumber
				\end{eqnarray}
								
				\begin{itemize}[label=$\bullet$, font=\normalsize, font=\color{black}, leftmargin=2cm,parsep=0cm, itemsep=0.25cm, topsep=0cm]
					\item[$ $] where $ i_{1} < i_{2}  ~...~ <i_{p} < r_{p} < r_{2} ~...~ < r_{1} \leq n $, $ r_{p} - i_{p} \geq 3 $ and $ p < n/2 $,
					\item[$ $] with  $ i_{p} < j $ and $ j +1  < r_{p}$,
					\item[$ $] while $d_{1} < d_{2} ~...~  < d_{z+1} $ and $ j + c \geq d_{ c+1} $, for $ 0 \leq c \leq z$.				
				\end{itemize}

				\begin{eqnarray}
					w_{2} &=& u_{1} a_{n+1} \sigma_{i_{1}} .. \sigma_{2} \sigma_{1} \sigma_{r_{1}} .. \sigma_{n-1} \sigma_{n} a_{n+1} \sigma_{i_{2}} .. \sigma_{2} \sigma_{1} \sigma_{r_{2}} .. \sigma_{n-1} \sigma_{n} a_{n+1} ~...~ \sigma_{i_{p}} .. \sigma_{2} \sigma_{1} \sigma_{r_{p}} .. \sigma_{n-1} \sigma_{n} \nonumber\\
					& & a_{n+1} \sigma_{j } .. \sigma_{2} \sigma_{1} \sigma_{j +2} .. \sigma_{n-1} \sigma_{n} \nonumber\\
					& & (a_{n+1} \sigma_{j +1} .. \sigma_{2} \sigma_{1}  \sigma_{j + 2} ..\sigma_{n-1} \sigma_{n})^{ m-(p+2)}\nonumber\\
					& & a_{n+1} \sigma_{j+1} ..\sigma_{d_{1}} \sigma_{j +2 } .. \sigma_{d_{2}} \sigma_{j +3} .. \sigma_{d_{3}} ~...~ \sigma_{j + z} ..\sigma_{d_{z}} \nonumber
				\end{eqnarray}
				
				\begin{itemize}[label=$\bullet$, font=\normalsize, font=\color{black}, leftmargin=2cm,parsep=0cm, itemsep=0.25cm, topsep=0cm]
					\item[$ $] where $  i_{1} < i_{2}  ~...~ <i_{p} < r_{p} < r_{2} ~...~ r_{1} \leq n $, $ r_{p} - i_{p} \geq 4 $, and $ p < n/2 $,
					\item[$ $] with $ i_{p} < j $ and $ j + 2 < r_{p}$,
					\item[$ $] while $ d_{1} < d_{2} ~...~  < d_{z} $ and $ j + c \geq d_{ c+1} $ for $ 0 \leq c \leq z $.				
				\end{itemize}
				
				\begin{eqnarray}
					w_{3} &=& u_{1} a_{n+1} \sigma_{i_{1}} .. \sigma_{2} \sigma_{1} \sigma_{r_{1}} .. \sigma_{n-1} \sigma_{n} a_{n+1} \sigma_{i_{2}} .. \sigma_{2} \sigma_{1} \sigma_{r_{2}} .. \sigma_{n-1} \sigma_{n} a_{n+1} ~...~ \sigma_{i_{p}} .. \sigma_{2} \sigma_{1} \sigma_{r_{p}}  \nonumber\\
					& & a_{n+1} \sigma_{l_{1}} .. \sigma_{g_{1}} \sigma_{l_{2}} .. \sigma_{g_{2}} ~...~ \sigma_{l_{t}} .. \sigma_{g_{t}}. \nonumber
				\end{eqnarray}
								
				\begin{itemize}[label=$\bullet$, font=\normalsize, font=\color{black}, leftmargin=2cm,parsep=0cm, itemsep=0.25cm, topsep=0cm]
					\item[$ $] where $  i_{1} < i_{2}  ~...~ <i_{p} < r_{p} < r_{2} ~...~ r_{1} \leq n $, $ r_{p} - i_{p} \geq 3 $ and $ p < n/2$,
					\item[$ $] with $ 1\leq l_{1} < l_{2} .. < l_{t} \leq n $ and $ 1\leq g_{1} < g_{2} .. < g_{t} \leq n $, 
					\item[$ $] while  $ i_{p} < l_{1} $, $g_{t} <r_{p} $ and $ g_{k} \leq l_{k}$ for any $1 \leq k \leq n $.\\
				\end{itemize}
			\end{lemme}

			\begin{demo}
				Before starting with the details of the proof, we call the reader's attention to the fact that our assumption that $\sigma_{n}$ is on the right in $ u_{h} $ for $ 2 \leq h \leq m$ is legitimate, since we know that these $ u_{h} $ are full by the discussion above. Using the discussion above we can write:
				\vspace{-0.25cm}
				\begin{eqnarray}
					u_{h-1} = \sigma_{i_{h}} ..\sigma_{2} \sigma_{1} \sigma_{r_{h}} ~...~\sigma_{n-1} \sigma_{n} x_{h}, \text{ for }3 \leq h \leq m+1. \nonumber
				\end{eqnarray}
				
				As above $1 \leq i_{h} \leq n-1 $, $ 1 \leq r_{h} \leq n $ and $i_{h} < r_{h}$, with $Supp(x_{h})\subseteq\{\sigma_{2},\sigma_{3}..\sigma_{n-1},\}$. Since $a_{n+1}$ commutes with $x_{h}$ for all $h$, we can write $x_{i}a_{n+1}u_{i+1}$ as $a_{n+1} u'_{i+1}$ with $u'_{i+1}$ full, in which $\sigma_{n}$ is on the right. Applying this inductively, we can write $w$ as follows: 
				\begin{eqnarray}
					w&=&u_{1} a_{n+1} \sigma_{i_{1}} .. \sigma_{2} \sigma_{1} \sigma_{r_{1}} .. \sigma_{n-1} \sigma_{n} a_{n+1} \sigma_{i_{2}} .. \sigma_{2} \sigma_{1} \sigma_{r_{2}} .. \sigma_{n-1} \sigma_{n} a_{n+1} ~...~ \sigma_{i_{m-1}} .. \sigma_{2} \sigma_{1} \sigma_{r_{m-1}} ~...~ \nonumber\\
					& & ~...~ \sigma_{n-1} \sigma_{n} a_{n+1} u _{m+1},\nonumber
				\end{eqnarray}
			
				with $u_{1}, u_{m+1}, i_{h}$ and $r_{h}$ as above. Now we have 3 main cases to consider: \\

				\begin{itemize}[label=$\bullet$, font=\normalsize, font=\color{black}, leftmargin=2cm,parsep=0cm, itemsep=0.25cm, topsep=0cm]
				
					\item[(1)] $r_{1} - i_{1} = 1 $, i.e., $r_{1} = i_{1} +1 $. \\ 
						
						In this case we do not have many choices for the full elements on the right of $u_{2}$: we have one and only one choice, $i_{h} = i_{1}$ for all $h \leq m-1$. Thus $j=i_{1}$. We have:
						\vspace{-0.25cm}
						\begin{eqnarray}
							w=u_{1} (a_{n+1}  \sigma_{j} ..\sigma_{2} \sigma_{1}  \sigma_{j +1} ..\sigma_{n-1} \sigma_{n})^{m-1} a_{n+1} u_{m+1}. \nonumber
						\end{eqnarray}
												
						Here we see that $u_{m+1}$ is a fully commutative element, which need not to be full, yet this element cannot have a reduced expression starting by any simple reflection in $W(A_{n})$ but  $\sigma_{i_{1}}$. If $u_{m+1}\neq 1$, we can thus, express it as follows:
						\begin{eqnarray}
							u_{m+1} &=& \sigma_{j} ..\sigma_{d_{1}} \sigma_{j+1} .. \sigma_{d_{2}} \sigma_{j +2} .. \sigma_{d_{3}} ~...~ \sigma_{j+z} .. \sigma_{d_{z+1}}, \nonumber
						\end{eqnarray}

						\vspace{-0.5cm}
						\begin{eqnarray}						
							\text{where }  d_{1} < d_{2} ~...~  < d_{z+1} \text{ and } j + c \geq d_{c+1} \text{ for } 0 \leq c \leq z.  \nonumber\\\nonumber
						\end{eqnarray}

					\item[(2)] $r_{1} - i_{1} = 2 $, i.e., $r_{1} = i_{1} +2 $.  \\ 
	
						In this case we have, as well, only one choice for the full element on the right of $u_{2}$, namely (we set $ i_{1} = j $):
						\begin{eqnarray}
							w&=&u_{1} a_{n+1} \sigma_{j} .. \sigma_{2} \sigma_{1} \sigma_{j +2} .. \sigma_{n-1} \sigma_{n} (a_{n+1} \sigma_{j +1} .. \sigma_{2} \sigma_{1} \sigma_{j + 2} .. \sigma_{n-1} \sigma_{n})^{ m-2} a_{n+1} u_{ m+1}, \nonumber
						\end{eqnarray}
						
						with conditions on $u_{m+1}$ analogous to those of case (1), that is:
						\begin{eqnarray}
							u_{m+1} &=& \sigma_{j+1} .. \sigma_{d_{1}} \sigma_{j+2} .. \sigma_{d_{2}} \sigma_{i_{1} +3} .. \sigma_{d_{3}} .. \sigma_{j + z} .. \sigma_{d_{z}}, \nonumber
						\end{eqnarray}	
						
						\vspace{-0.5cm}
						\begin{eqnarray}
							\text{where }  d_{1} < d_{2} ~...~  < d_{z} \text{ and } j + c \geq d_{c} \text{ for } 1 \leq c \leq z.  \nonumber\\\nonumber
						\end{eqnarray}
						
					\item[(3)] $ r_{1} - i_{1} > 2 $. 
						\vspace{0.15cm}
						\begin{eqnarray}
							\text{Say } w&=&u_{1} a_{n+1} \sigma_{i_{1}} .. \sigma_{2} \sigma_{1} \sigma_{r_{1}} .. \sigma_{n-1} \sigma_{n} a_{n+1} \sigma_{i_{2}} .. \sigma_{2} \sigma_{1} \sigma_{r_{2}} .. \sigma_{n-1} \sigma_{n}\nonumber\\
							& & a_{n+1} ~...~ \sigma_{i_{m-1}} .. \sigma_{2} \sigma_{1} \sigma_{r_{m-1}} ~...~ \sigma_{n-1} \sigma_{n} a_{n+1} u _{m+1},\nonumber\\\nonumber
						\end{eqnarray}
				\end{itemize}
	
				we see that we have to choose $ r_{2} $ and $i_{2}$ such that $i_{1} < i_{2} < r_{2}  < r_{1} $. Hence, after a finite  number of steps, we will face one of the cases  (1) or (2). Thus we have one of the next forms:\\ 
				
				\begin{itemize}[label=$\bullet$, font=\normalsize, font=\color{black}, leftmargin=2cm,parsep=0cm, itemsep=0.25cm, topsep=0cm]
	
					\item[(1')] This is the case related to (1), i.e., we have:
						\begin{eqnarray}
							w &=& u_{1} a_{n+1} \sigma_{i_{1}} .. \sigma_{2} \sigma_{1} \sigma_{r_{1}} .. \sigma_{n-1} \sigma_{n} a_{n+1} \sigma_{i_{2}} .. \sigma_{2} \sigma_{1} \sigma_{r_{2}} .. \sigma_{n-1} \sigma_{n} a_{n+1} ~...~ \sigma_{i_{p}} .. \sigma_{2} \sigma_{1} \sigma_{r_{p}} .. \sigma_{n-1} \sigma_{n}\nonumber\\
							& & (a_{n+1} \sigma_{j} ..\sigma_{2} \sigma_{1} \sigma_{j +1} ..\sigma_{n-1} \sigma_{n})^{m-(1+p)} a_{n+1} u_{m+1}, \nonumber
						\end{eqnarray}
					
						here $ i_{1} < i_{2}  ~...~ <i_{p} < r_{p} < r_{2} ~...~< r_{1} \leq n $ and $ r_{p} - i_{p} \geq 3 $. We have necessarily $ p < n/2 $, while $u_{m+1}$ is as in case (1). \\
					
					\item[(2')] This case is related to (2), i.e., we have:		
						\begin{eqnarray}
							w&=&u_{1} a_{n+1} \sigma_{i_{1}} .. \sigma_{2} \sigma_{1} \sigma_{r_{1}} .. \sigma_{n-1} \sigma_{n} a_{n+1} \sigma_{i_{2}} .. \sigma_{2} \sigma_{1} \sigma_{r_{2}} .. \sigma_{n-1} \sigma_{n} a_{n+1} ~...~ \sigma_{i_{p}} .. \sigma_{2} \sigma_{1} \sigma_{r_{p}} .. \sigma_{n-1} \sigma_{n}\nonumber\\
							& & a_{n+1} \sigma_{j } .. \sigma_{2} \sigma_{1} \sigma_{j +2} .. \sigma_{n-1} \sigma_{n} (a_{n+1}	\sigma_{j +1} .. \sigma_{2} \sigma_{1} \sigma_{j + 2} .. \sigma_{n-1} \sigma_{n})^{ m-(p+2)} a_{n+1} u_{ m+1}, \nonumber
						\end{eqnarray}
						
						here $ i_{1} < i_{2}  ~...~ <i_{p} < r_{p} < r_{2} ~...~< r_{1} \leq n $ and $ r_{p} - i_{p} \geq 4 $. We have necessarily $ p < n/2 $, while $u_{m+1}$ is as in case (2).  \\  
					\vspace{0.25cm}	

					\item[(3')] This case is related to some "short" elements (with respect to $L$): \\
					
						suppose that we stopped picking pairs $ (i , r) $ before having a difference of $1$ or $2$ between them, hence:
						\vspace{-0.25cm}
						\begin{eqnarray}
							w&=&u_{1} a_{n+1} \sigma_{i_{1}} .. \sigma_{2} \sigma_{1} \sigma_{r_{1}} .. \sigma_{n-1} \sigma_{n} a_{n+1} \sigma_{i_{2}} .. \sigma_{2} \sigma_{1} \sigma_{r_{2}} .. \sigma_{n-1} \sigma_{n} a_{n+1} ~...~ \sigma_{i_{p}} .. \sigma_{2} \sigma_{1} \sigma_{r_{p}}	a_{n+1} u_{ m+1},\nonumber
						\end{eqnarray}

						with  $ i_{1} < i_{2}  ~...~ <i_{p} < r_{p} < r_{2} ~...~ <r_{1} \leq n $ and $ r_{p} - i_{p} \geq 3 $. We have necessarily $ p < n/2 $. \\ 
						
				\end{itemize}	
				
				In this case, the choice of $u_{m+1}$ is much more complicated than in the other two cases. It has the form: 
				\vspace{-0.25cm}
				\begin{eqnarray}
					\sigma_{l_{1}} .. \sigma_{g_{1}} \sigma_{l_{2}} .. \sigma_{g_{2}} ~...~ \sigma_{l_{t}} .. \sigma_{g_{t}}, \nonumber
				\end{eqnarray}
				
		       \end{demo}
				
					 where $ 1\leq l_{1} < l_{2} .. < l_{t} \leq n $, $ 1\leq g_{1} < g_{2} .. < g_{t} \leq n $ and $ g_{k} \leq l_{k}$ for any $1 \leq k \leq n $.  And in addition we have $ i_{p} < l_{1} $ and $g_{t} <r_{p}$. \\

				 Later on, we will be back to handle the possible forms of $u_{m+1}$ in details.

			\begin{definition}
				In elements of type $w_{1}$, the following element is called the short block:
				\begin{eqnarray}
					a_{n+1} \sigma_{i_{1}} .. \sigma_{2} \sigma_{1} \sigma_{r_{1}} .. \sigma_{n-1} \sigma_{n} a_{n+1} \sigma_{i_{2}} .. \sigma_{2} \sigma_{1} \sigma_{r_{2}} .. \sigma_{n-1} \sigma_{n} a_{n+1} ~...~ \sigma_{i_{p}} .. \sigma_{2} \sigma_{1} \sigma_{r_{p}} .. \sigma_{n-1} \sigma_{n}. \nonumber
				\end{eqnarray}
				
				We call $(a_{n+1} \sigma_{j} .. \sigma_{2} \sigma_{1} \sigma_{j+1} .. \sigma_{n-1} \sigma_{n})^{m-(1+p)}$ the convergent block of $w_{1}$.\\
	
				We call $a_{n+1} u_{m+1} $  the residue block of $w_{1}$.\\
	
				Hence we can write $w_{1} = u_{1}$. short block. convergent block. residue block. (We do the same thing for elements of type $w_{2}$, in which for example, the convergent block is $(a_{n+1} \sigma_{j +1} .. \sigma_{2}\sigma_{1} \sigma_{j + 2} .. \sigma_{n-1} \sigma_{n})^{ m-(p+2)}$ ).\\		
			\end{definition}			
					
			\begin{definition}
				An element of the last two types is called short, if and only if its convergent block is equal to 1.\\
			\end{definition}
	
			\begin{remarque}\label{2_4_14} 
				It is easy to see that $ w_{1} $ and $w_{2} $ could be unified in the following form:
				\begin{eqnarray}
					w_{1} &=& u_{1} a_{n+1} \sigma_{i_{1}} .. \sigma_{2} \sigma_{1} \sigma_{r_{1}} .. \sigma_{n-1} \sigma_{n} a_{n+1} \sigma_{i_{2}} .. \sigma_{2} \sigma_{1} \sigma_{r_{2}} .. \sigma_{n-1} \sigma_{n} a_{n+1} ~...~ \sigma_{i_{p}} .. \sigma_{2} \sigma_{1} \sigma_{r_{p}} .. \sigma_{n-1} \sigma_{n}\nonumber\\
					& & \big(a_{n+1} \sigma_{j} .. \sigma_{2} \sigma_{1} \sigma_{j+1} .. \sigma_{n-1} \sigma_{n}\big)^{m-\big(1+L(\text{the~short block})\big)} \nonumber\\
					& & a_{n+1} \sigma_{j} .. \sigma_{d_{1}} \sigma_{j +1 } .. \sigma_{d_{2}} \sigma_{i_{1}+2} .. \sigma_{d_{3}} ~...~ \sigma_{i_{1}+z} .. \sigma_{d_{z+1}}, \nonumber
				\end{eqnarray}
								
				\begin{itemize}[label=$\bullet$, font=\normalsize, font=\color{black}, leftmargin=2cm,parsep=0cm, itemsep=0.25cm, topsep=0cm]
					\item[$ $] where: $ i_{1} < i_{2} ~...~ <i_{p} < r_{p} < r_{2} ~...~ < r_{1} \leq n $, $ r_{p} - i_{p} \geq 2 $, and $ p $ necessarily lesser than $n/2 $. 
					\item[$ $] while $ i_{p} < j $, $ j +1  < r_{p} $, $ d_{1} < d_{2} ~...~  < d_{z+1}$ and $ j + c \geq d_{ c+1} $ for $ 0 \leq c \leq z $.\\
					
				\end{itemize}
				
			Nevertheless, for the moment, we will go on keeping looking at them as two different forms .\\
			
			\end{remarque}		
	
			We see that the set of short elements is of finite cardinal, because  of the fact that the affine length $L$ of such elements is bounded.  Special cases of the last lemma, which comes from the 3 types above when $ u_{m+1} = 1 $, are included in the general formula.\\
			
			Now we classify the elements of $W^{c}(\tilde{A_{n}})$ with $ n\geq 3 $.\\ 
	 
			Consider an arbitrary $w$ in  $W^{c}(\tilde{A_{n}})$ with $L(w) \geq 2 $, written as:
			\begin{eqnarray}
				w = u_{1} a_{n+1} u_{2} a_{n+1} ~...~ u_{m} a_{n+1} u_{m+1}.  \nonumber
			\end{eqnarray}
			
			We start the classification, depending on the choice of $u_{1}$ which can have one, and only one of the following forms:\\
		
		 \begin{itemize}[label=$\bullet$, font=\normalsize, font=\color{black}, leftmargin=2cm,parsep=0cm, itemsep=0.25cm, topsep=0cm]          
			\item[$(a)$] $ u_{1}$ is full, with  $ \sigma_{n} $ on the left.
			\item[$(b)$] $ u_{1}$ is full, with  $ \sigma_{n} $ on the left. 
			\item[$(c)$] $ \sigma_{n} $ belongs to $ supp (u_{1})$  and  $ \sigma_{1} $ does not.
			\item[$(d)$] $ \sigma_{1} $ belongs to $ supp (u_{1})$  and  $ \sigma_{n} $ does not.
			\item[$(e)$] $u_{1} = 1 $.\\
		\end{itemize} 
	        
			\clearpage
			Suppose that we are in case $(a)$. \\
			
			We have $u_{1} =  \sigma_{n} \sigma_{n-1} ~...~ \sigma_{2}\sigma_{1}$. In this case there is only one choice for the full elements $u_{i}$  with $ 2 \leq i \leq m$, which is to be equal to $u_{1}$, hence $w = (\sigma_{n} \sigma_{n-1} ~...~ \sigma_{2}\sigma_{1} a_{n+1})^{m} u_{m+1} $. Here $u_{m+1}$ is either 1 or $\sigma_{n} \sigma_{n-1} ~...~ \sigma_{j} $, thus we have two possible types:
			\begin{eqnarray}
				x_{1} &=& (\sigma_{n} \sigma_{n-1} ~...~ \sigma_{2}\sigma_{1} a_{n+1})^{m} \sigma_{n} \sigma_{n-1} ~...~ \sigma_{i}, \text{ for } 1 \leq 
		i \leq n. \nonumber\\
				x_{2} &=& (\sigma_{n} \sigma_{n-1} ~...~ \sigma_{2}\sigma_{1} a_{n+1})^{m} .\nonumber			
			\end{eqnarray}
		
			Suppose that we are in case (b).\\
			
			Set $u_{1} := \sigma_{i_{0}} ..\sigma_{2} \sigma_{1} \sigma_{r_{0}} ~...~ \sigma_{n-1} \sigma_{n} x_{0}$.  It is clear that $ u_{i }$, for $ 2\leq i$, cannot be equal to $\sigma_{n} ~...~ \sigma_{1}$, hence all the full elements $u_{i}$, for $2 \leq i \leq m$, have $\sigma_{n}$ on the right. Here we can use the same discussion as in lemma.2.4.9. We arrive to the possible types, by replacing $ u_{1} a_{n+1}$ (in which $w$ starts) by $1$. Thus we have three possible types (modulo maybe a shift of indexes to the left):
			\begin{eqnarray}
				x_{3} &=& \sigma_{i_{1}} .. \sigma_{2} \sigma_{1} \sigma_{r_{1}} .. \sigma_{n-1} \sigma_{n} a_{n+1} \sigma_{i_{2}} .. \sigma_{2} \sigma_{1} \sigma_{r_{2}} .. \sigma_{n-1} \sigma_{n} a_{n+1} ~...~ \sigma_{i_{p}} .. \sigma_{2} \sigma_{1} \sigma_{r_{p}} .. \sigma_{n-1} \sigma_{n} \nonumber\\
				& & (a_{n+1} \sigma_{j} .. \sigma_{2} \sigma_{1} \sigma_{j +1 } .. \sigma_{n-1} \sigma_{n})^{m-p}  \nonumber\\
				& & a_{n+1} \sigma_{j} .. \sigma_{d_{1}} \sigma_{j+1} .. \sigma_{d_{2}} \sigma_{j+2} .. \sigma_{d_{3}} ~...~ \sigma_{j+z} .. \sigma_{d_{z+1}},\nonumber
			\end{eqnarray}
								
			\begin{itemize}[label=$\bullet$, font=\normalsize, font=\color{black}, leftmargin=2cm,parsep=0cm, itemsep=0.25cm, topsep=0cm]
				\item[$ $] where  $ i_{1} < i_{2} ~...~ <i_{p} < r_{p} < r_{2} ~...~< r_{1} \leq n $, $ r_{p} - i_{p} \geq 3 $ and $ p < n/2 $,
				\item[$ $] with $ i_{p} < j $ and $ j +1  < r_{p} $,
				\item[$ $] while  $ d_{1} < d_{2} ~...~  < d_{z+1}$ and $ j + c \geq d_{ c+1} $ for $ 0 \leq c \leq z $. 	
			\end{itemize}		
						
			\begin{eqnarray}
				x_{4} &=& \sigma_{i_{1}} .. \sigma_{2} \sigma_{1} \sigma_{r_{1}} .. \sigma_{n-1} \sigma_{n} a_{n+1} \sigma_{i_{2}} .. \sigma_{2} \sigma_{1} \sigma_{r_{2}} .. \sigma_{n-1} \sigma_{n} a_{n+1} ~...~ \sigma_{i_{p}} .. \sigma_{2} \sigma_{1} \sigma_{r_{p}} .. \sigma_{n-1} \sigma_{n}~~~~\nonumber\\
				& & a_{n+1} \sigma_{j} .. \sigma_{2} \sigma_{1} \sigma_{j +2} .. \sigma_{n-1} \sigma_{n} (a_{n+1} \sigma_{j+1} .. \sigma_{2} \sigma_{1} \sigma_{j + 2} .. \sigma_{n-1} \sigma_{n})^{ m-(p+1)}  \nonumber\\
				& & a_{n+1} \sigma_{j+1} .. \sigma_{d_{1}} \sigma_{j +2 } .. \sigma_{d_{2}} \sigma_{j+3} .. \sigma_{d_{3}} ~...~ \sigma_{j+z} .. \sigma_{d_{z}},\nonumber
			\end{eqnarray}
			\vspace{-0.35cm}
								
			\begin{itemize}[label=$\bullet$, font=\normalsize, font=\color{black}, leftmargin=2cm,parsep=0cm, itemsep=0.25cm, topsep=0cm]
				\item[$ $] where  $ i_{1} < i_{2} ~...~ <i_{p} < r_{p} < r_{2} ~...~< r_{1} \leq n $, $ r_{p} - i_{p} \geq 4  $ and $ p < n/2 $,
				\item[$ $] with $ i_{p} < j $ and $ j +2  < r_{p} $,
				\item[$ $] while $ d_{1} < d_{2} ~...~  < d_{z} $ and $ i_{1} + c \geq d_{ c} $ for $ 1 \leq c \leq z$.
			\end{itemize}		

			\begin{eqnarray}
				x_{5} &=& \sigma_{i_{1}} .. \sigma_{2} \sigma_{1} \sigma_{r_{1}} .. \sigma_{n-1} \sigma_{n} a_{n+1} \sigma_{i_{2}} .. \sigma_{2} \sigma_{1} \sigma_{r_{2}} .. \sigma_{n-1} \sigma_{n} a_{n+1} ~...~ \sigma_{i_{p}} .. \sigma_{2} \sigma_{1} \sigma_{r_{p}} ~~~~~~~~~~~~~~~~\nonumber\\
			 & & a_{n+1} \sigma_{l_{1}} .. \sigma_{g_{1}} \sigma_{l_{2}} .. \sigma_{g_{2}} ~...~ \sigma_{l_{t}} .. \sigma_{g_{t}},\nonumber
			\end{eqnarray}
			
			\vspace{0.1cm}
			\begin{itemize}[label=$\bullet$, font=\normalsize, font=\color{black}, leftmargin=2cm,parsep=0cm, itemsep=0.25cm, topsep=0cm]
				\item[$ $]where $ i_{1} < i_{2}  ~...~ <i_{p} < r_{p} < r_{2} ~...~ < r_{1} \leq n $, $ r_{p} - i_{p} \geq 3 $ and $ p <n/2$,
				\item[$ $] with $ 1\leq l_{1} < l_{2} .. < l_{t} \leq n $ and $ 1\leq g_{1} < g_{2} .. < g_{t} \leq n $,
				\item[$ $] while $ i_{p} < l_{1} $, $g_{t} <r_{p} $ and $ g_{k} \leq l_{k}$ for any $1 \leq k \leq n $. \\ 
	
			\end{itemize}
			 
			Of course, we keep in mind the three special cases $ x'_{3}$, (resp. $x'_{4}$ and $ x'_{5} $), which are obtained from $ x_{3}$, (resp. $x_{4}$ and $ x_{5} $) by replacing $u_{m+1}$ by $1$.\\ 
			
			Suppose that we are in case (c). \\
		
			Here, $u_{1}$ can be written as $ u_{1} = \sigma_{h}\sigma_{h+1} .. \sigma_{n-1}\sigma_{n} y $, where $ 2 \leq h \leq n $, $y$ is in $W (A_{n-1}) $ and $ \sigma_{1} \notin supp (y) $. Hence, we can write $ w $ as follows:
			\begin{eqnarray}
				w = \sigma_{h} .. \sigma_{n} a_{n+1} u_{2} a_{n+1} ~...~ u_{m} a_{n+1} u_{m+1}. \nonumber
			\end{eqnarray}
	
			We see that $u_{2} $ (thus every $ u_{i} $ with $ 2 \leq i \leq m+1$) cannot start with $\sigma_{n}$. That means each $ u_{i} $, with $ 2 \leq i \leq m $, is a full element in which $\sigma_{n} $ is on the left. By using the lemma \ref{2_4_11} $w$ is one of the three following elements:
			\begin{eqnarray}
				x_{6} &=& \sigma_{h} .. \sigma_{n} a_{n+1} \sigma_{i_{1}} ..\sigma_{2} \sigma_{1} \sigma_{r_{1}} ..\sigma_{n-1} \sigma_{n} a_{n+1} \sigma_{i_{2}} ..\sigma_{2} \sigma_{1}  \sigma_{r_{2}} ..\sigma_{n-1} \sigma_{n}  a_{n+1} ~...~ \sigma_{i_{p}} ..\sigma_{2} \sigma_{1} \sigma_{r_{p}} ~...~\sigma_{n-1} \sigma_{n} \nonumber\\
				& & (a_{n+1} \sigma_{j} ..\sigma_{2} \sigma_{1} \sigma_{j +1 } ..\sigma_{n-1} \sigma_{n})^{m-(1+p)} \nonumber\\
				& & a_{n+1} \sigma_{j} ..\sigma_{d_{1}} \sigma_{j +1 } .. \sigma_{d_{2}} \sigma_{j +2} .. \sigma_{d_{3}} ~...~ \sigma_{j+ z} .. \sigma_{d_{z+1}}, \nonumber			
			\end{eqnarray}
			\vspace{0.1cm}
								
			\begin{itemize}[label=$\bullet$, font=\normalsize, font=\color{black}, leftmargin=2cm,parsep=0cm, itemsep=0.25cm, topsep=0cm]
				\item[$ $] where $i_{1} < i_{2} ~...~ <i_{p} < r_{p} < r_{2} ~...~< r_{1} \leq n $, $ r_{p} - i_{p} \geq 3 $ and $ p < n/2 $,
				\item[$ $] with $ i_{p} < j $, $ j +1 < r_{p}$, $ d_{1} < d_{2} ~...~  < d_{z+1} $  and $ j + c \geq d_{ c+1} $ for $ 0 \leq c \leq z$,
				\item[$ $] while $ i_{1} < h  $, and if $ r_{1} - i_{1} > 1$, then $ h < r_{1}$.
			\end{itemize}	
									
			\begin{eqnarray}
				x_{7} &=& \sigma_{h} .. \sigma_{n} a_{n+1} \sigma_{i_{1}} ..\sigma_{2} \sigma_{1} \sigma_{r_{1}} ..\sigma_{n-1} \sigma_{n} a_{n+1} \sigma_{i_{2}} .. \sigma_{2} \sigma_{1} \sigma_{r_{2}} ..\sigma_{n-1} \sigma_{n}  a_{n+1} ~...~ \sigma_{i_{p}} ..\sigma_{2} \sigma_{1} \sigma_{r_{p}} ~...~\sigma_{n-1} \sigma_{n} \nonumber\\
				& & a_{n+1} \sigma_{i_{1}} ..\sigma_{2} \sigma_{1} \sigma_{i_{1}+2} ..\sigma_{n-1} \sigma_{n} \nonumber\\
				& & (a_{n+1} \sigma_{j +1} .. \sigma_{2} \sigma_{1} \sigma_{j + 2} ..\sigma_{n-1} \sigma_{n})^{ m-(p+2)} \nonumber\\
				& & a_{n+1} \sigma_{j+1} ..\sigma_{d_{1}} \sigma_{j +2 } .. \sigma_{d_{2}} \sigma_{j +3} .. \sigma_{d_{3}} ~...~ \sigma_{j+ z} ..\sigma_{d_{z}}, \nonumber			
			\end{eqnarray}
								
			\begin{itemize}[label=$\bullet$, font=\normalsize, font=\color{black}, leftmargin=2cm,parsep=0cm, itemsep=0.25cm, topsep=0cm]
				\item[$ $] where $i_{1} < i_{2} ~...~ <i_{p} < r_{p} < r_{2} ~...~ < r_{1} \leq n $, $ r_{p} - i_{p} \geq 4 $ and $ p < n/2 $,
				\item[$ $] with $ i_{p} < i_{1} $, $i_{1} + 2 < r_{p}$, $ d_{1} < d_{2} ~...~ < d_{z} $ and $ i_{1} + c \geq d_{ c}$ for $ 1 \leq c \leq z$,
				\item[$ $] while $ i_{1} < h  $, and if $ r_{1} - i_{1} > 1$  then $ h < r_{1}$.
			\end{itemize}
			
			\begin{eqnarray}
				x_{8} &=& \sigma_{h} .. \sigma_{n} a_{n+1} \sigma_{i_{1}} ..\sigma_{2} \sigma_{1} \sigma_{r_{1}} ..\sigma_{n-1} \sigma_{n} a_{n+1} \sigma_{i_{2}} ..\sigma_{2} \sigma_{1} \sigma_{r_{2}} ..\sigma_{n-1} \sigma_{n}  a_{n+1} ~...~ \sigma_{i_{p}} ..\sigma_{2} \sigma_{1} \sigma_{r_{p}} ~~~~~~~~\nonumber\\
				& & a_{n+1} \sigma_{l_{1}} .. \sigma_{g_{1}} \sigma_{l_{2}} .. \sigma_{g_{2}} ~...~ \sigma_{l_{t}} .. \sigma_{g_{t}}, \nonumber			
			\end{eqnarray}
								
			\begin{itemize}[label=$\bullet$, font=\normalsize, font=\color{black}, leftmargin=2cm,parsep=0cm, itemsep=0.25cm, topsep=0cm]
				\item[$ $] here $ i_{1} < i_{2}  ~...~ <i_{p} < r_{p} < r_{2} ~...~ < r_{1} \leq n $, $ r_{p} - i_{p} \geq 3 $ and $ p < n/2 $,
				\item[$ $] where $ 1\leq l_{1} < l_{2} .. < l_{t} \leq n $ and $ 1\leq g_{1} < g_{2} .. < g_{t} \leq n $,
				\item[$ $] with $ g_{k} \leq l_{k}$ for any $1 \leq k \leq n $, $ i_{p} < l_{1} $ and $g_{t} <r_{p} $,
				\item[$ $] while $ i_{1} < h  $, and if $ r_{1} - i_{1} > 1$  then $ h < r_{1}$.\\
			\end{itemize}	
				
			As before we keep in mind the three special cases $ x'_{6}$, (resp. $x'_{7}$ and $ x'_{8} $), which are obtained from $ x_{6}$, (resp. $x_{7}$ and $ x_{8} $) by replacing $u_{m+1}$ by $1$.\\
	
			Suppose that we are in case (d). \\	

			Here, $ u_{1} $ can be written $ \sigma_{h} \sigma_{h-1} .. \sigma_{1} y $, where $y$ is in $W (A_{n-1}) $, with $ \sigma_{1} \notin supp (y) $ and $1 \leq h \leq n-1 $. Hence we can suppose that
			\begin{eqnarray}
				w = \sigma_{h} .. \sigma_{1} a_{n+1} u_{2} a_{n+1} ~...~ u_{m} a_{n+1} u_{m+1}. \nonumber		
			\end{eqnarray}
			
			Here, we have two main choices for $ u_{2} $. The first one is that $ \sigma_{n} $ is on the left, then $w$ has the following form:
			\begin{eqnarray}
				x_{9} &=& \sigma_{h} .. \sigma_{1} a_{n+1} (\sigma_{n} \sigma_{n-1} ~...~ \sigma_{2}\sigma_{1} a_{n+1})^{m} \sigma_{n}  ~...~ \sigma_{i}, \text{ where } 1 \leq  h \leq n-1 \text{, and } 1 \leq  i \leq n. \nonumber			
			\end{eqnarray}

			The second choice is that $u_{i}$, for $ 2 \leq i \leq n $, has $\sigma_{n}$ on the right. As above we have three forms, namely:
			\begin{eqnarray}
				x_{10} &=& \sigma_{h} .. \sigma_{1} a_{n+1} \sigma_{i_{1}} ..\sigma_{2} \sigma_{1} \sigma_{r_{1}} ..\sigma_{n-1} \sigma_{n}a_{n+1} \sigma_{i_{2}} ..\sigma_{2} \sigma_{1} \sigma_{r_{2}} ..\sigma_{n-1} \sigma_{n}  a_{n+1} ~...~ \sigma_{i_{p}} ..\sigma_{2} \sigma_{1} \sigma_{r_{p}} ~...~\sigma_{n-1} \sigma_{n}\nonumber \\
				& & (a_{n+1} \sigma_{j} ..\sigma_{2} \sigma_{1} \sigma_{j +1 } ..\sigma_{n-1} \sigma_{n})^{m-(1+p)} \nonumber \\
				& & a_{n+1} \sigma_{j} ..\sigma_{d_{1}} \sigma_{j +1 } .. \sigma_{d_{2}} \sigma_{j +2} .. \sigma_{d_{3}} ~...~ \sigma_{j + z} ..\sigma_{d_{z+1}}, \nonumber	
			\end{eqnarray}
		
			\begin{itemize}[label=$\bullet$, font=\normalsize, font=\color{black}, leftmargin=2cm,parsep=0cm, itemsep=0.25cm, topsep=0cm]
				\item[$ $] where $  i_{1} < i_{2} ~...~ <i_{p} < r_{p} < r_{2} ~...~ < r_{1} \leq n $ and $ r_{p} - i_{p} \geq 3 $ and $ p < n/2$,
				\item[$ $] with $i_{p} < i_{1}$, $ i_{1} +1 < r_{p} $, $d_{1} < d_{2} ~...~  < d_{z+1}$ and $i_{1} + c \geq d_{ c+1}$ for $0 \leq c \leq z$,
				\item[$ $] while $ i_{1} < h  $, and if $ r_{1} - i_{1} > 1$,  then $ h < r_{1}$.
			\end{itemize}
			
			\vspace{-0.05cm}
			\begin{eqnarray}
				x_{11} &=& \sigma_{h} .. \sigma_{1} a_{n+1} \sigma_{i_{1}} ..\sigma_{2} \sigma_{1} \sigma_{r_{1}} ..\sigma_{n-1} \sigma_{n} a_{n+1} \sigma_{i_{2}} .. \sigma_{2} \sigma_{1} \sigma_{r_{2}} ..\sigma_{n-1} \sigma_{n}  a_{n+1} ~...~ \sigma_{i_{p}} ..\sigma_{2} \sigma_{1} \sigma_{r_{p}} ~...~\sigma_{n-1} \sigma_{n} \nonumber\\
				& & a_{n+1} \sigma_{i_{1}} ..\sigma_{2} \sigma_{1} \sigma_{i_{1} +2 } ..\sigma_{n-1} \sigma_{n} \nonumber\\
				& &(a_{n+1} \sigma_{j +1} .. \sigma_{2} \sigma_{1} \sigma_{j + 2} ..\sigma_{n-1} \sigma_{n})^{ m-(p+2)} \nonumber\\
				& &a_{n+1} \sigma_{j+1} ..\sigma_{d_{1}} \sigma_{j +2 } .. \sigma_{d_{2}} \sigma_{j +3} .. \sigma_{d_{3}} ~...~ \sigma_{j + z} ..	\sigma_{d_{z}}, \nonumber
			\end{eqnarray}			
		
			\begin{itemize}[label=$\bullet$, font=\normalsize, font=\color{black}, leftmargin=2cm,parsep=0cm, itemsep=0.25cm, topsep=0cm]
				\item[$ $] where $ i_{1} < i_{2}  ~...~ <i_{p} < r_{p} < r_{2} ~...~ < r_{1} \leq n $, $ r_{p} - i_{p} \geq 4 $ and $ p < n/2 $,
				\item[$ $] with $i_{p} < i_{1} $, $ i_{1} + 2 < r_{p} $, $d_{1} < d_{2} ~...~  < d_{z}$ and $i_{1} + c \geq d_{ c}$ for $ 1 \leq c \leq z$, 
				\item[$ $] while $ i_{1} < h  $, and if $ r_{1} - i_{1} > 1$,  then $ h < r_{1}$.
			\end{itemize}	
			
			\begin{eqnarray}
				x_{12} &=& \sigma_{h} .. \sigma_{1} a_{n+1} \sigma_{i_{1}} ..\sigma_{2} \sigma_{1} \sigma_{r_{1}} ..\sigma_{n-1} \sigma_{n} a_{n+1} \sigma_{i_{2}} ..\sigma_{2} \sigma_{1} \sigma_{r_{2}} ..\sigma_{n-1} \sigma_{n}  a_{n+1} ~...~ \sigma_{i_{p}} ..\sigma_{2} \sigma_{1}  \sigma_{r_{p}}~~~~~~~~~~~~~~\nonumber\\ 
				& & a_{n+1} \sigma_{l_{1}} .. \sigma_{g_{1}} \sigma_{l_{2}} .. \sigma_{g_{2}} ~...~ \sigma_{l_{t}} .. \sigma_{g_{t}}, \nonumber
			\end{eqnarray}
		
			\begin{itemize}[label=$\bullet$, font=\normalsize, font=\color{black}, leftmargin=2cm,parsep=0cm, itemsep=0.25cm, topsep=0cm]
				\item[$ $] where $ i_{1} < i_{2}  ~...~ <i_{p} < r_{p} < r_{2} ~...~ < r_{1} \leq n $, $ r_{p} - i_{p} \geq 3 $ and $ p < n/2 $,
				\item[$ $] with $ 1\leq l_{1} < l_{2} .. < l_{t} \leq n $ and $ 1\leq g_{1} < g_{2} .. < g_{t} \leq n $,
				\item[$ $] while $ g_{k} \leq l_{k}$, for any $1 \leq k \leq n $, with $ i_{p} < l_{1} $, $g_{t} <r_{p} $ and $ i_{1} < h < r_{1} $.\\
			\end{itemize}
			
			Still, we keep in mind the three special cases $ x'_{9}$, (resp. $x'_{10}$ and $ x'_{11} $), which are obtained from $ x_{9}$, (resp. $x_{10}$ and $ x_{11} $) by replacing $u_{m+1}$ by $1$.\\
			
			Suppose that we are in case (e). \\	
			
			This case will be a particular case of the above cases. We use the following notation in  $W(\tilde{A_{n}})$: $ \sigma_{0} = \sigma_{n+1} = 1$. With this notation we see that types $ x_{1}, x_{2} $ and $ x_{9} $ could be unified in one form, say $ c_{1}$.\\

			Moreover, $x_{3}$ (resp. $x_{4}$ and $x_{5}$) can be unified in one form with $x_{6}$ (resp. $x_{7}$ and $x_{8}$), when $ i_{1} = 0 $.\\ 
			
			Similarly, $x_{3}$ (resp. $x_{4}$ and $x_{5}$) can be unified in one form with $x_{10}$ (resp. $x_{11}$ and $x_{12}$), when $ r_{1} = n+1 $.\\
							
			From what precedes, we formulate our classification by the following corollary.\\
			
			\begin{corollaire}
				Let $n \geq 3 $. Let $w$ be in $W^{c}(\tilde{A_{n}})$, such that $ 2 \leq L (w) $. Then $ w $ has one of the following forms:	
				\begin{eqnarray}
					c_{1} &=& \sigma_{j} ~...~ \sigma_{2} \sigma_{1} (\sigma_{n} \sigma_{n-1} ~...~ \sigma_{2}\sigma_{1} a_{n+1})^{m} \sigma_{n} \sigma_{n-1} ~...~ \sigma_{i}, \nonumber\\
					& & \text{ where } 1 \leq i \leq n+1 \text{ and } 0 \leq j \leq n .\nonumber\\\nonumber\\
					c_{2} &=& \sigma_{i_{1}} ..\sigma_{2} \sigma_{1} \sigma_{r_{1}} ..\sigma_{n-1} \sigma_{n}  a_{n+1} \sigma_{i_{2}} .. \sigma_{2} \sigma_{1} \sigma_{r_{2}} ..\sigma_{n-1} \sigma_{n} a_{n+1} ~...~ \sigma_{i_{p}} ..\sigma_{2} \sigma_{1} \sigma_{r_{p}} ~...~\sigma_{n-1} \sigma_{n}\nonumber\\
					& & (a_{n+1} \sigma_{j} ..\sigma_{2} \sigma_{1} \sigma_{j +1 } ..\sigma_{n-1} \sigma_{n})^{m-(p)} \nonumber\\
					& & a_{n+1} \sigma_{j} ..\sigma_{d_{1}} \sigma_{j +1 } .. \sigma_{d_{2}} \sigma_{j +2} .. \sigma_{d_{3}} ~...~ \sigma_{j + z} ..\sigma_{d_{z+1}},\nonumber
				\end{eqnarray} 
				\vspace{-0.25cm}
		
				\begin{itemize}[label=$\bullet$, font=\normalsize, font=\color{black}, leftmargin=2cm,parsep=0cm, itemsep=0.25cm, topsep=0cm]
					\item[$ $] where  $ 0 \leq i_{1} < i_{2}~...~ <i_{p} < r_{p} < r_{2} ~...~< r_{1} \leq n +1  $ and $ r_{p} - i_{p} \geq 3$, 
					\item[$ $] with $ i_{p} < j $ and $ j +1  < r_{p} $,
					\item[$ $] while $ d_{1} < d_{2} ~...~  < d_{z+1} $ and $ j + c \geq d_{ c+1} $ for  $ 0 \leq c \leq z $.
				\end{itemize}
				
				\begin{eqnarray}
					c_{3} &=& \sigma_{i_{1}} ..\sigma_{2} \sigma_{1} \sigma_{r_{1}} ..\sigma_{n-1} \sigma_{n}  a_{n+1} \sigma_{i_{2}} .. \sigma_{2} \sigma_{1} \sigma_{r_{2}} ..\sigma_{n-1} \sigma_{n} a_{n+1} ~...~ \sigma_{i_{p}} ..\sigma_{2} \sigma_{1} \sigma_{r_{p}} ~...~\sigma_{n-1} \sigma_{n} \nonumber\\
					& & a_{n+1} \sigma_{j } ..\sigma_{2} \sigma_{1} \sigma_{j +2 } ..\sigma_{n-1} \sigma_{n} \nonumber\\
					& & (a_{n+1} \sigma_{j +1} .. \sigma_{2} \sigma_{1} \sigma_{j + 2} ..\sigma_{n-1} \sigma_{n})^{ m-(p+1)} \nonumber\\
					& & a_{n+1} \sigma_{j+1} ..\sigma_{d_{1}} \sigma_{j +2 } .. \sigma_{d_{2}} \sigma_{j +3} .. \sigma_{d_{3}} ~...~ \sigma_{j + z} ..\sigma_{d_{z}}, \nonumber							
				\end{eqnarray}		
				\vspace{-0.25cm}
		
				\begin{itemize}[label=$\bullet$, font=\normalsize, font=\color{black}, leftmargin=2cm,parsep=0cm, itemsep=0.25cm, topsep=0cm]
					\item[$ $] where  $ 0 \leq i_{1} < i_{2}  ~...~ <i_{p} < r_{p} < r_{2} ~...~ < r_{1} \leq n+1 $ and $ r_{p} - i_{p} \geq 4 $,
					\item[$ $] with $ i_{p} < j $ and $ j + 2 < r_{p} $,
					\item[$ $] while $ d_{1} < d_{2} ~...~  < d_{z} $ and $ i_{1} + c \geq d_{ c} $ for $ 1 \leq c \leq z$.  					
				\end{itemize}	
	
				\begin{eqnarray}
					c_{4} &=& \sigma_{i_{1}} ..\sigma_{2} \sigma_{1} \sigma_{r_{1}} ..\sigma_{n-1} \sigma_{n}  a_{n+1} \sigma_{i_{2}} ..\sigma_{2} \sigma_{1} \sigma_{r_{2}} ..\sigma_{n-1} \sigma_{n} a_{n+1} ~...~ \sigma_{i_{p}} ..\sigma_{2} \sigma_{1}  \sigma_{r_{p}} ~~~~~~~~~~~~~~\nonumber\\
					& & a_{n+1} \sigma_{l_{1}} .. \sigma_{g_{1}} \sigma_{l_{2}} .. \sigma_{g_{2}} ~...~ \sigma_{l_{t}} .. \sigma_{g_{t}}, \nonumber
				\end{eqnarray}
				\vspace{-0.25cm}
		
				\begin{itemize}[label=$\bullet$, font=\normalsize, font=\color{black}, leftmargin=2cm,parsep=0cm, itemsep=0.25cm, topsep=0cm]
					\item[$ $] where: $ 0 \leq i_{1} < i_{2}  ~...~ <i_{p} < r_{p} < r_{2} ~...~ < r_{1} \leq n $ and $ r_{p} - i_{p} \geq 3 $, 
					\item[$ $] with $ 1\leq l_{1} < l_{2} .. < l_{t} \leq n $ and $ 1\leq g_{1} < g_{2} .. < g_{t} \leq n $, 
					\item[$ $] while $ i_{p} < l_{1} $, $g_{t} <r_{p} $ and $ g_{k} \leq l_{k}$ for any $1 \leq k \leq n $.\\
				\end{itemize}
	
				In all cases $p$ is necessarily bounded by $ n/2 $ .	
			\end{corollaire}
			
			In order to get to the final form of our classification we shall do one more step, explained in the following remark:
			\vspace{0.25cm}
			\begin{remarques}
				We set $ \sigma_{t} = 1 $ when $ 0 \geq t$ or $t \geq n+1  $. We can actually unify cases $ c_{1}, c_{2} $ and $ c_{3} $ with the following formula:
				\begin{eqnarray}
					& & \sigma_{i_{1}} ..\sigma_{2} \sigma_{1} \sigma_{r_{1}} ..\sigma_{n-1} \sigma_{n} a_{n+1} \sigma_{i_{2}} .. \sigma_{2} \sigma_{1} \sigma_{r_{2}} ..\sigma_{n-1} \sigma_{n} a_{n+1} ~...~ \sigma_{i_{p}} ..\sigma_{2} \sigma_{1} \sigma_{r_{p}} ~...~\sigma_{n-1} \sigma_{n}\nonumber\\
					& & (a_{n+1} \sigma_{j} ..\sigma_{2} \sigma_{1} \sigma_{j +1 } ..\sigma_{n-1} \sigma_{n})^{K} \nonumber\\
					& & a_{n+1} \sigma_{j} ..\sigma_{d_{1}} \sigma_{j +1 } .. \sigma_{d_{2}} \sigma_{j +2} .. \sigma_{d_{3}} ~...~ \sigma_{j + z} ..\sigma_{d_{z+1}},\nonumber
				\end{eqnarray}
				\vspace{-0.5cm}
		
				\begin{itemize}[label=$\bullet$, font=\normalsize, font=\color{black}, leftmargin=2cm,parsep=0cm, itemsep=0.25cm, topsep=0cm]
					\item[$ $] where $ 0 \leq i_{1} < i_{2}~...~ <i_{p} < r_{p} < r_{2} ~...~< r_{1} \leq n +1  $ and $ r_{p} - i_{p} \geq 2 $,
					\item[$ $] with $ i_{p} < j $, $ j  \leq r_{p} -1 $, $i_{1} \leq n $ and $ 1\leq K $,
					\item[$ $] while  $ d_{1} < d_{2} ~...~  < d_{z+1} $ and $ j + c \geq d_{ c+1} $ for $ 0 \leq c \leq z $. \\
				\end{itemize}		
			\end{remarques}

			Moreover, we see that our formula expresses the elements of $W^{c}(\tilde{A_{2}}) $. With this last remark, the proof of theorem \ref{2_4_4} is completed, after noticing that the way in which we get the general form, ensures the unicity of this form. \\ 
			
	\section{Fully commutative affine braids} \label{2_5}
	    \vspace{0.3cm}

		Now we consider the tower of affine braid groups:\\
		\begin{eqnarray}
			B(\tilde{A_{0}}) \stackrel{F_{0}}{\longrightarrow} B(\tilde{A_{1}})\stackrel{F_{1}}{\longrightarrow} ~...~ B(\tilde{A_{n-1}}) \stackrel{F_{n}}{\longrightarrow}B(\tilde{A_{n}}) \stackrel{F_{n+1}}{\longrightarrow}~...~ \nonumber\\\nonumber
		\end{eqnarray}
		
  		where $B(\tilde{A_{0}})$ is the trivial group. Via $F_{n}$, every $ B(\tilde{A_{n-1}})$ injects into  $B(\tilde{A_{n}})$, where $F_{n}$ is induced by the injection of the $B$-type braid groups $B(B_{n}) \hookrightarrow B(B_{n+1})$, noticing that $B(\tilde{A_{n}})$ is a subgroup of $B(B_{n+1})$ for $n \geq 0$ (see \cite{Graham_Lehrer_2003}). The injection  $F_{n}$ is given as follows: 
		
		       \begin{eqnarray}
					F_{n}: B(\tilde{A_{n-1}}) &\longrightarrow& B(\tilde{A_{n}})  \nonumber\\
					\sigma_{i} &\longmapsto& \sigma_{i}$ ~~~ \text{for} $1\leq i\leq n-1 \nonumber\\
					a_{n} &\longmapsto& \sigma_{n} a_{n+1}\sigma^{-1}_{n}. \nonumber\\\nonumber
				\end{eqnarray}
		
		We are interested with viewing $ B(\tilde{A_{n}})$ containing $ B(\tilde{A_{n-1}})$. The following computations are done in view of understanding the tower of affine Temperley-Lieb algebras.\\

  		 In what follows we give a kind of normal form for the lift, in $B(\tilde{A_{n}})$, of fully commutative elements in $ W(\tilde{A_{n}}) $ see corollary \ref{2_5_2}. We keep the same symbols for the generators of the affine braid group and their images via the natural surjection onto affine Coxeter group.\\
	
		Let $1 \leq n$. Let $\bar{w}$ be in $ W^{c}(\tilde{A_{n}})$. The general form of $\bar{w}$ is
		
		\begin{eqnarray}
			\bar{w} &=& \sigma_{i_{1}} ..\sigma_{2} \sigma_{1} \sigma_{r_{1}} ..\sigma_{n-1} \sigma_{n} a_{n+1} \sigma_{i_{2}} .. \sigma_{2} \sigma_{1} \sigma_{r_{2}} ..\sigma_{n-1} \sigma_{n} a_{n+1} ~...~ \sigma_{i_{p}} ..\sigma_{2} \sigma_{1} \sigma_{r_{p}} ~...~\sigma_{n-1} \sigma_{n} \nonumber\\
			& & (a_{n+1} \sigma_{j} ..\sigma_{2} \sigma_{1} \sigma_{j +1 } ..\sigma_{n-1} \sigma_{n})^{k} \bar{u}, \nonumber
		\end{eqnarray}	
		\vspace{-0.5cm}
		
		\begin{itemize}[label=$\bullet$, font=\normalsize, font=\color{black}, leftmargin=2cm,parsep=0cm, itemsep=0.25cm, topsep=0cm]
			\item[$ $] where $ 0 \leq i_{1} < i_{2}~...~ <i_{p} < r_{p} < r_{2} ~...~< r_{1} \leq n +1  $ and $ r_{p} - i_{p} \geq 2, $
			\item[$ $] with $ i_{p} < j $, $ j  \leq r_{p} -1 $, $i_{1} \leq n $ and $ 0 \leq k $,
			\item[$ $] while $\bar{u} = a_{n+1} \bar{v} $, with $v$ fully commutative in  $ W({A_{n}}) $.\\
		\end{itemize}
	
	 	We lift $\bar{w}$,(resp. $\bar{u}$ and $ \bar{v}$) to $ w$,(resp. $u$ and $v$) in $B(\tilde{A_{n}})$. Assume that $j<n$, i.e., $w$ is not of the form $v(a_{n+1} \sigma_{n} .. \sigma_{2} \sigma_{1})^{k} u $. We show that $w$ has the form:\\
	 	
	 	$h(\sigma_{n} .. \sigma_{1} a_{n+1})^{m} x$, where $x$ is in $B(A_{n})$ and $h$ is in $ B(\tilde{A_{n-1}}) $. \\
		
		We show as well that $\bar{w} $ has the form:\\
		
		$ \bar{f}(\sigma_{n} .. \sigma_{1} a_{n+1})^{m} \sigma_{n} .. \sigma_{i}$, where $1 \leq i \leq n+1$, $m$ is a positive integer and $\bar{f}$ is in $W(\tilde{A_{n-1}})$.\\
	
		Recall that $ \sigma_{n} \sigma_{n-1}.. \sigma_{1}a_{n+1} $ acts on the elements of $ B(\tilde{A_{n-1}}) $ exactly the way $\phi^{-1}_{n}$ does in $B(B_{n})$. In order to simplify, we set $\phi^{-1}_{n}= \psi$. We write $(\sigma_{n} .. \sigma_{1} a_{n+1})^{d}h = \psi^{d} \left[ h\right]  (\sigma_{n} .. \sigma_{1} a_{n+1})^{d} $, for any $h$ in $ B(\tilde{A_{n-1}}) $. The automorphism  $\phi$ is of order $n$, of course, $d$ is to be taken mod$n$. We keep in mind that $a_{n+1} a_{n}= a_{n} \sigma_{n}= \sigma_{n}a_{n+1} $.

		\begin{lemme} \label{2_5_1}
			Let $y$ be $ \sigma_{j} ..\sigma_{2} \sigma_{1} \sigma_{j +1 } ..\sigma_{n-1} \sigma_{n} a_{n+1}$ in $ B(\tilde{A_{n}}) $. Let $w= y^{k}$ with $2 \leq j \leq n-1$ and $ 2\leq k $. Suppose that $k= m(n-j+1) +r $ where $0\leq r < n-j+1 $. Then: \\
			
			\clearpage
			
			\begin{itemize}[label=$\bullet$, font=\normalsize, font=\color{black}, leftmargin=2cm,parsep=0cm, itemsep=0.25cm, topsep=0cm]
				\item[$(1)$] If $ m=0 $ we have:
					\begin{eqnarray}
						w = (\sigma_{j} ..\sigma_{2} \sigma_{1} \sigma_{j +1 } ..\sigma_{n-1}  a_{n})^{r}\sigma_{n} \sigma_{n-1} .. \sigma_{n+1-r}. \nonumber
					\end{eqnarray}
					\vspace{-0.5cm}

				\item[$(2)$] if $ 0<m $ we have:
					\begin{eqnarray}
						w &=& \prod^{i=m-1}_{i=0} \psi^{i} \left[ (\sigma_{j} ..\sigma_{2} \sigma_{1} \sigma_{j +1} ..\sigma_{n-1} a_{n})^{n-j}\sigma_{j}.. \sigma_{n-1} \right] \psi^{m}\big[(\sigma_{j} ..\sigma_{2} \sigma_{1}  \sigma_{j +1 } ..\sigma_{n-1} a_{n})^{r}\big]\nonumber\\
						& & (\sigma_{n} \sigma_{n-1} ..\sigma_{2} \sigma_{1} a_{n+1})^{m} \sigma_{n} \sigma_{n-1} .. \sigma_{n+1-r}.\nonumber		
					\end{eqnarray}			
			\end{itemize}
		\end{lemme}
				
		\begin{demo}
			
			\begin{eqnarray}
				\text{We have }y &=& \sigma_{j} ..\sigma_{2} \sigma_{1} \sigma_{j +1 } ..\sigma_{n-1} \sigma_{n}a_{n+1} =  \sigma_{j} ..\sigma_{2} \sigma_{1} \sigma_{j +1 } ..\sigma_{n-1} a_{n} \sigma_{n}. \nonumber \\
				\text{Hence, }y^{2} &=& \sigma_{j} ..\sigma_{2} \sigma_{1} \sigma_{j +1 } ..\sigma_{n-1} a_{n}\sigma_{n}  \sigma_{j} ..\sigma_{2} \sigma_{1} \sigma_{j +1 } ..\sigma_{n-1} \sigma_{n}a_{n+1} \nonumber\\ 
				&=& \sigma_{j} ..\sigma_{2} \sigma_{1} \sigma_{j +1 } ..\sigma_{n-1} a_{n} \sigma_{j} ..\sigma_{2} \sigma_{1}	\sigma_{j +1 } ..\sigma_{n-1} \sigma_{n}a_{n+1}\sigma_{n-1} \nonumber\\
				&=& (\sigma_{j} ..\sigma_{2} \sigma_{1} \sigma_{j +1 } ..\sigma_{n-1} a_{n})^{2}\sigma_{n} \sigma_{n-1}. \nonumber
    			\end{eqnarray}
		\vspace{-0.5cm}
			Continuing this way, we can see that  whenever $ 0 \leq r \leq n-j $:\\
			\begin{eqnarray}
				y^{r} &=& (\sigma_{j} ..\sigma_{2} \sigma_{1} \sigma_{j +1 } ..\sigma_{n-1} a_{n})^{r}\sigma_{n} \sigma_{n-1} .. \sigma_{s}, \text{ with } s+r = n+1. \text{ in particular:} \nonumber\\ \nonumber\\
				y^{n-j} &=& (\sigma_{j} ..\sigma_{2} \sigma_{1} \sigma_{j +1 } ..\sigma_{n-1} a_{n})^{n-j}\sigma_{n} \sigma_{n-1}.. \sigma_{j+1}.\text{ Thus:} \nonumber\\\nonumber\\
				y^{n-j+1} &=& (\sigma_{j} ..\sigma_{2} \sigma_{1} \sigma_{j +1 } ..\sigma_{n-1}  a_{n})^{n-j}\sigma_{n} \sigma_{n-1}.. \sigma_{j+1}\sigma_{j} ..\sigma_{2} \sigma_{1} \sigma_{j +1 } ..\sigma_{n-1} \sigma_{n} a_{n+1}  \nonumber\\\nonumber\\
				&=& (\sigma_{j} ..\sigma_{2} \sigma_{1} \sigma_{j +1 } ..\sigma_{n-1} a_{n})^{n-j}\sigma_{n} \sigma_{n-1}.. \sigma_{j+1}\sigma_{j} ..\sigma_{2} \sigma_{1}\sigma_{j +1 } .. \sigma_{n-1}\sigma_{n}  a_{n+1}\nonumber\\\nonumber\\
				&=& (\sigma_{j} ..\sigma_{2} \sigma_{1} \sigma_{j +1 } ..\sigma_{n-1} a_{n})^{n-j}\sigma_{j} .. \sigma_{n-1} \sigma_{n} \sigma_{n-1} ..\sigma_{2} \sigma_{1} a_{n+1} \text{. We see that: } \nonumber\\\nonumber\\ 
				y^{2(n-j+1)} &=& \psi^{0}\left[(\sigma_{j} ..\sigma_{2} \sigma_{1} \sigma_{j +1 } ..\sigma_{n-1}  a_{n})^{n-j} \sigma_{j } .. \sigma_{n-1}\right] \nonumber\\\nonumber\\
				& & .\psi^{1}\left[(\sigma_{j} ..\sigma_{2} \sigma_{1} \sigma_{j +1 } ..\sigma_{n-1} F(a_{n}))^{n-j}\sigma_{j}.. \sigma_{n-1}\right] (\sigma_{n} \sigma_{n-1} ..\sigma_{2} \sigma_{1} a_{n+1})^{2}.\nonumber
			\end{eqnarray}
	
			In the same way, for $m\geq 0$, considering the action of $ \sigma_{n} \sigma_{n-1} ..\sigma_{2} \sigma_{1} a_{n+1}$ on $B(\tilde{A_{n-1}})$, we see that: 
			\begin{eqnarray}
				y^{m(n-j+1)} = \prod^{i=m-1}_{i=0} \psi^{i} \left[ (\sigma_{j} ..\sigma_{2} \sigma_{1} \sigma_{j +1 } ..\sigma_{n-1} F(a_{n}))^{n-j} \sigma_{j}.. \sigma_{n-1}\right] (\sigma_{n} \sigma_{n-1} ..\sigma_{2} \sigma_{1} a_{n+1})^{m} .\nonumber
			\end{eqnarray}
			
			Finally, let $k= m(n-j+1) +r $, where $0\leq r < n-j+1 $. We have:
			\begin{eqnarray}
				y^{k} &=& \prod^{i=m-1}_{i=0} \psi^{i} \left[ (\sigma_{j} ..\sigma_{2} \sigma_{1} \sigma_{j +1 } ..\sigma_{n-1} F(a_{n}))^{n-j}\sigma_{j}.. \sigma_{n-1} \right]\nonumber\\\nonumber\\
				& & (\sigma_{n} \sigma_{n-1} ..\sigma_{2} \sigma_{1} a_{n+1})^{m} (\sigma_{j} ..\sigma_{2} \sigma_{1} \sigma_{j +1 } ..\sigma_{n-1}  a_{n})^{r}\sigma_{n} \sigma_{n-1} .. \sigma_{n+1-r}\text{. Thus,} \nonumber\\\nonumber\\
				y^{k} &=& \prod^{i=m-1}_{i=0} \psi^{i} \left[ (\sigma_{j} ..\sigma_{2} \sigma_{1} \sigma_{j +1 } ..\sigma_{n-1} F(a_{n}))^{n-j}\sigma_{j}.. \sigma_{n-1} \right] \psi^{m}\left[(\sigma_{j} ..\sigma_{2} \sigma_{1} \sigma_{j +1 } ..\sigma_{n-1}  a_{n})^{r}\right]\nonumber\\\nonumber\\
				& & (\sigma_{n} \sigma_{n-1} ..\sigma_{2} \sigma_{1} a_{n+1})^{m} \sigma_{n} \sigma_{n-1} .. \sigma_{n+1-r}.\nonumber	
			\end{eqnarray}
				
		\end{demo}
		
		In particular, for $j=1$, i.e., $ w= (  \sigma_{1} \sigma_{2} .. \sigma_{n} a_{n+1})^{k} $, we have:
		\begin{eqnarray}
			y^{k} = \prod^{i=m-1}_{i=0} \psi^{i} \left[ (\sigma_{1} .. \sigma_{n-1} a_{n})^{n-1} \sigma_{1} .. \sigma_{n-1}\right] \psi^{m} \left[ (\sigma_{1} .. \sigma_{n-1} a_{n})^{r} \right] (\sigma_{n} .. \sigma_{1} a_{n+1})^{m} \sigma_{n} \sigma_{n-1} .. \sigma_{i}. \nonumber
		\end{eqnarray}
	
		Now we go back to the general form of $\bar{w}$, that is: 
		\begin{eqnarray}
			\bar{w} &=& \sigma_{i_{1}} ..\sigma_{2} \sigma_{1} \sigma_{r_{1}} ..\sigma_{n-1} \sigma_{n} a_{n+1} \sigma_{i_{2}} .. \sigma_{2} \sigma_{1} \sigma_{r_{2}} ..\sigma_{n-1} \sigma_{n} a_{n+1} ~...~ \sigma_{i_{p}} ..\sigma_{2} \sigma_{1} \sigma_{r_{p}} ~...~\sigma_{n-1} \sigma_{n} \nonumber\\
			& & (a_{n+1} \sigma_{j} ..\sigma_{2} \sigma_{1} \sigma_{j +1 } ..\sigma_{n-1} \sigma_{n})^{k} \bar{u}, \nonumber
		\end{eqnarray}	
		\vspace{-0.5cm}
		
		\begin{itemize}[label=$\bullet$, font=\normalsize, font=\color{black}, leftmargin=2cm,parsep=0cm, itemsep=0.25cm, topsep=0cm]
			\item[$ $] where $ 0 \leq i_{1} < i_{2}~...~ <i_{p} < r_{p} < r_{2} ~...~< r_{1} \leq n +1  $ and $ r_{p} - i_{p} \geq 2, $
			\item[$ $] with $ i_{p} < j $, $ j  \leq r_{p} -1 $, $i_{1} \leq n $ and $ 0 \leq k $,
			\item[$ $] while $\bar{u} = a_{n+1} \bar{v} $, with $v$ fully commutative in  $ W({A_{n}}) $.\\
		\end{itemize}
		
		\vspace{-1cm}
		
		\begin{eqnarray}
			\text{Then, } w &=& \sigma_{i_{1}} ..\sigma_{2} \sigma_{1}  \sigma_{r_{1}} ..\sigma_{n-1} \sigma_{n}  a_{n+1} \sigma_{i_{2}} .. \sigma_{2} \sigma_{1} \sigma_{r_{2}} ..\sigma_{n-1} \sigma_{n} a_{n+1} ~...~ \sigma_{i_{p}} ..\sigma_{2} \sigma_{1} \sigma_{r_{p}} ~...~\sigma_{n-1} \sigma_{n} a_{n+1} \nonumber\\
			& & (\sigma_{j} ..\sigma_{2} \sigma_{1} \sigma_{j +1 } ..\sigma_{n-1} \sigma_{n} a_{n+1})^{k}  v, \nonumber
		\end{eqnarray}
		\vspace{-0.4cm}
		
		\begin{itemize}[label=$\bullet$, font=\normalsize, font=\color{black}, leftmargin=2cm,parsep=0cm, itemsep=0.25cm, topsep=0cm]
			
			\item[$ $] with conditions similar to those above. In particular, $v$ is fully commutative in  $ W({A_{n}})$.\\
		\end{itemize}
		
		We see that $ r_{i} \leq r_{1} - i $, but $ r_{1} \leq n+1 $. which gives: $  r_{i} \leq n- i+1 $ for all $i$.\\ 
	
		Here we treat two main cases:\\
			
		\begin{itemize}[label=$\bullet$, font=\normalsize, font=\color{black}, leftmargin=2cm,parsep=0cm, itemsep=0.25cm, topsep=0cm]
			\item $ r_{1} \leq n $, that is $\sigma_{n}$ belongs to the support of $\sigma_{i_{1}} ..\sigma_{2} \sigma_{1} \sigma_{r_{1}} ..\sigma_{n-1} \sigma_{n} $.
			\item $ r_{1} = n+1 $, that is $\sigma_{n}$ belongs to the support of $\sigma_{i_{1}} ..\sigma_{2} \sigma_{1}  \sigma_{r_{1}} ..\sigma_{n-1} \sigma_{n} $, (this case covers the case where $p=0$).\\
		\end{itemize} 
				
		We start by the first case $ r_{1} \leq n $ (we suppose that $ p> 0$).	
		\begin{eqnarray}
			\text{Set } x := \sigma_{i_{1}} ..\sigma_{2} \sigma_{1} \sigma_{r_{1}} ..\sigma_{n-1} \sigma_{n} a_{n+1} \sigma_{i_{2}} .. \sigma_{2} \sigma_{1} \sigma_{r_{2}} ..\sigma_{n-1} \sigma_{n} a_{n+1} ~...~ \sigma_{i_{p}} ..\sigma_{2} \sigma_{1}  \sigma_{r_{p}} ~...~\sigma_{n-1} \sigma_{n} a_{n+1}.  \nonumber
		\end{eqnarray}
	
		Thus, $w = x(\sigma_{j} ..\sigma_{2} \sigma_{1}  \sigma_{j +1 } ..\sigma_{n-1} \sigma_{n} a_{n+1})^{k}v$. Repeating the first step  of lemma 2.5.1 (keeping in mind that $ r_{i} \leq n-i $), we see that:
		\begin{eqnarray}
			x &=& \sigma_{i_{1}} ..\sigma_{2} \sigma_{1} \sigma_{r_{1}} ..\sigma_{n-1} a_{n} \sigma_{i_{2}} .. \sigma_{2} \sigma_{1} \sigma_{r_{2}} ..\sigma_{n-1} \sigma_{n} a_{n+1} ~...~ \sigma_{i_{p}} ..\sigma_{2} \sigma_{1} \sigma_{r_{p}} ~...~\sigma_{n-1} \sigma_{n} a_{n+1} \sigma_{n-(p-1)} \nonumber\\ 
			 &=& \sigma_{i_{1}} ..\sigma_{2} \sigma_{1} \sigma_{r_{1}} ..\sigma_{n-1} a_{n} \sigma_{i_{2}} ..	\sigma_{2} \sigma_{1}  \sigma_{r_{2}} ..\sigma_{n-1} a_{n} ~...~ \sigma_{i_{p}} ..\sigma_{2} \sigma_{1}  \sigma_{r_{p}} ~...~\sigma_{n-1} \sigma_{n} a_{n+1} \sigma_{n-p} \sigma_{n-(p-1)}.\nonumber
		\end{eqnarray}
		
		After $p$ steps we see that:
		\begin{eqnarray}
			x = \sigma_{i_{1}} ..\sigma_{2} \sigma_{1}  \sigma_{r_{1}} ..\sigma_{n-1} a_{n}  \sigma_{i_{2}} ..\sigma_{2} \sigma_{1}  \sigma_{r_{2}} ..\sigma_{n-1} a_{n}  ..  \sigma_{i_{p}} ..\sigma_{2} \sigma_{1} \sigma_{r_{p}} ~...~\sigma_{n-1} a_{n} \sigma_{n}\sigma_{n-1} .. \sigma_{n - (p-1)}. \nonumber
		\end{eqnarray}
		
		Now, we set $ \epsilon:= n - (p-1) $. We show that $ \epsilon > j+1 $ as follows. \\ 
	
		We know that $ j +1  \leq r_{p}  $, but $ r_{p} < n-p < n-(p-1) $. In other words:
		\begin{eqnarray}
			j+1 \leq r_{p} < n-(p-1), \text{ that is } j+1 < \epsilon . \nonumber
		\end{eqnarray}		
		\vspace{-0.75cm}
		\begin{eqnarray}
			\text{Set } \rho &:=& \sigma_{i_{1}} ..\sigma_{2} \sigma_{1}  \sigma_{r_{1}} ..\sigma_{n-1} a_{n} \sigma_{i_{2}} ..\sigma_{2} \sigma_{1} \sigma_{r_{2}} ..\sigma_{n-1} a_{n}  ~...~ \sigma_{i_{p}} ..\sigma_{2} \sigma_{1} \sigma_{r_{p}} ~...~\sigma_{n-1} a_{n}.\text{ Thus, }\nonumber\\\nonumber\\
			 w &=& \rho \sigma_{n} \sigma_{n-1} ..  \sigma_{\epsilon}(\sigma_{j} ..\sigma_{2} \sigma_{1}  \sigma_{j +1} ..\sigma_{n-1}\sigma_{n} a_{n+1})^{k} v, \text{ with } \rho \in B(\tilde{A_{n-1}}) \text{ and } \epsilon > j+1. \nonumber
		\end{eqnarray}
	
		Now we have $ w = \rho  \sigma_{n} \sigma_{n-1} .. \sigma_{\epsilon}y^{k} v$. Every $y$ acts on $ \sigma_{i} $ in the following way: $ \sigma_{i} y = y \sigma_{i-1} $, for  $ \epsilon \leq i \leq n$, since $ j+1 < \epsilon $ and hence $j+1 <i$.\\ 
	
		If $ k =0 $, the job is done (this case is included in the general form). \\ 
	
		Let $ 1\leq k $. We have two main cases:\\ 
		
		\begin{itemize}[label=$\bullet$, font=\normalsize, font=\color{black}, leftmargin=2cm,parsep=0cm, itemsep=0.25cm, topsep=0cm]
			\item[$(1)$] $1 \leq k < \epsilon - (j+1)$.
			\item[$(2)$] $\epsilon - (j+1) \leq k $.\\ 
		\end{itemize}
			
		We start by (1). Set $ e:= \epsilon - k $. We have:
		\begin{eqnarray}
			\sigma_{n} \sigma_{n-1} .. \sigma_{\epsilon} y^{k} &=& y^{k} \sigma_{e+(n-\epsilon)} \sigma_{e+(n-\epsilon) -1} ..	\sigma_{e}. \text{ That is,} \nonumber\\ 
			\sigma_{n} \sigma_{n-1} .. \sigma_{\epsilon} y^{k} &=& y^{k} \sigma_{n-k} \sigma_{n-1-k} ..\sigma_{\epsilon - k }. \nonumber 
		\end{eqnarray}
	
	     We have $ k < \epsilon - (j+1) = n- (p-1) - (j+1) = n-j-p < n-j $. Thus, in the terms of lemma \ref{2_5_1}, we are in case (1), (even with the same $y$). That is: 
		\begin{eqnarray}
			y^{k}= (\sigma_{j} ..\sigma_{2} \sigma_{1}  \sigma_{j +1 } ..\sigma_{n-1}  a_{n})^{k}\sigma_{n} \sigma_{n-1} ..\sigma_{n+1-k}. \nonumber
		\end{eqnarray}

		\vspace{-0.25cm}
		\begin{eqnarray}
			\text{Thus }w = \rho (\sigma_{j} ..\sigma_{2} \sigma_{1}  \sigma_{j +1 } ..\sigma_{n-1}  a_{n})^{k}\sigma_{n} \sigma_{n-1} .. \sigma_{n-1-k}  \sigma_{n-k}  \sigma_{n-1-k} ..  \sigma_{\epsilon - k } v \nonumber
		\end{eqnarray}
	
		which is basically the case $m=0$.\\
	
		Case (2), where  $\epsilon - (j+1)  \leq k $.\\
	
		We see that for $k = \epsilon- (j+1) $, we have:
		\begin{eqnarray}
			& & \sigma_{n}\sigma_{n-1} .. \sigma_{\epsilon} y^{k} = y^{k} \sigma_{n-k} \sigma_{(n-k)-1} .. \sigma_{j+1}. \text{ In other terms:}\nonumber\\\nonumber\\
			& & \sigma_{n}\sigma_{n-1} .. \sigma_{\epsilon} (y)^{\epsilon- (j+1)} = (y)^{\epsilon- (j+1)} \sigma_{j+p}\sigma_{j+p-1} ..\sigma_{j+1}. \nonumber\\\nonumber\\
			& &\text{That is } \sigma_{n}\sigma_{n-1} .. \sigma_{\epsilon} y^{k} = \sigma_{n} \sigma_{n-1} .. \sigma_{\epsilon} y^{\epsilon - (j+1)} y^{k-(\epsilon - (j+1))}. \nonumber\\\nonumber
		\end{eqnarray}
		
		Now, $\epsilon - (j+1) = n - (p-1)-(j+1) = n-j-p < n-j $. Here, we can apply lemma \ref{2_5_1} (again we are in the first case). Precisely :
		\begin{eqnarray}
			(y)^{\epsilon- (j+1)} &=& (\sigma_{j} ..\sigma_{2} \sigma_{1} \sigma_{j +1 } ..\sigma_{n-1} a_{n})^{\epsilon- (j+1)}\sigma_{n} \sigma_{n-1} .. \sigma_{n+1-(\epsilon- (j+1))} \nonumber\\
			&=& (\sigma_{j} ..\sigma_{2} \sigma_{1} \sigma_{j +1 } ..\sigma_{n-1} a_{n})^{\epsilon- (j+1)}\sigma_{n} \sigma_{n-1} .. \sigma_{j+p+1}. \nonumber
		\end{eqnarray}
		
		Set $h:= k-(\epsilon - (j+1))$. We get: 
		\begin{eqnarray}
			\sigma_{n}\sigma_{n-1} .. \sigma_{\epsilon}y^{k} =&(& \sigma_{j} ..\sigma_{2} \sigma_{1} \sigma_{j +1} ..\sigma_{n-1} a_{n})^{\epsilon - (j+1)}\sigma_{n} \sigma_{n-1} .. \sigma_{j+p+1} \sigma_{j+p}..\sigma_{j+1}yy^{h-1}.\nonumber\\\nonumber\\
			\text{But,}& &\sigma_{n} \sigma_{n-1} .. \sigma_{j+p+1} \sigma_{j+p}..\sigma_{j+1}yy^{h-1} =	\sigma_{n} .. \sigma_{j+1} \sigma_{j} ..\sigma_{2} \sigma_{1} \sigma_{j +1} ..\sigma_{n}a_{n+1}y^{h-1}, \nonumber\\ \nonumber\\
			\text{which is equal to } & & \sigma_{n} \sigma_{n-1} .. \sigma_{1}\underbrace{\sigma_{j+1} .. \sigma_{n-1} \sigma_{n}}_{} a_{n+1} y^{h-1}, \nonumber\\\nonumber\\ 
			\text{Which is } & & \underbrace{\sigma_{j} .. \sigma_{n-2} \sigma_{n-1}}_{}\sigma_{n} \sigma_{n-1} .. \sigma_{1}a_{n+1} y^{h-1}. \nonumber
		\end{eqnarray}
	
		Now set $\eta:= \rho  \sigma_{j} .. \sigma_{n-2} \sigma_{n-1} \in B(\tilde{A_{n-1}})  $. We get $w = \eta \sigma_{n} \sigma_{n-1} .. \sigma_{1}a_{n+1} y^{h-1} v.\\$

		\begin{itemize}[label=$\bullet$, font=\normalsize, font=\color{black}, leftmargin=2cm,parsep=0cm, itemsep=0.25cm, topsep=0cm]
			\item[$(a)$] If $ h-1 \leq n-j $, we see that:
				\begin{eqnarray}
					w &=& \eta \sigma_{n} \sigma_{n-1} .. \sigma_{1}a_{n+1} (\sigma_{j} ..\sigma_{2} \sigma_{1} \sigma_{j +1 } ..\sigma_{n-1} a_{n})^{h-1}\sigma_{n} \sigma_{n-1} .. \sigma_{n+1-(h-1)} v \nonumber \\
					 &=& \eta \psi \left[ (\sigma_{j} ..\sigma_{2} \sigma_{1} \sigma_{j +1 } ..\sigma_{n-1} a_{n})^{h-1} \right] \sigma_{n} \sigma_{n-1} .. \sigma_{1}a_{n+1} \sigma_{n} \sigma_{n-1} .. \sigma_{n+1-(h-1)}v. \nonumber\\	\nonumber
				\end{eqnarray}
			\item[$(b)$] If $n-j< h-1 $, we see that:
				\begin{eqnarray}
					w &=& \eta \sigma_{n} \sigma_{n-1} .. \sigma_{1}a_{n+1} \prod^{i=m-1}_{i=0} \psi^{i} \left[ (\sigma_{j} ..\sigma_{2} \sigma_{1} \sigma_{j +1 } ..\sigma_{n-1}  a_{n})^{n-j}\sigma_{j}.. \sigma_{n-1} \right]\nonumber\\
					& & \psi^{m}\left[(\sigma_{j} ..\sigma_{2} \sigma_{1} \sigma_{j +1} ..\sigma_{n-1} a_{n})^{r}\right] (\sigma_{n} \sigma_{n-1} ..\sigma_{2} \sigma_{1} a_{n+1})^{m} \sigma_{n} \sigma_{n-1} .. \sigma_{n+1-r} v .\nonumber\\\nonumber\\
					 \text{Thus, }w&=& \eta  \prod^{i=m}_{m=1} \psi^{i} \left[ (\sigma_{j} ..\sigma_{2}\sigma_{1} \sigma_{j +1 } ..\sigma_{n-1}  a_{n})^{n-j}\sigma_{j}.. \sigma_{n-1} \right] \psi^{m+1}\left[(\sigma_{j} ..\sigma_{2} \sigma_{1} \sigma_{j +1 } ..\sigma_{n-1} a_{n})^{r}\right] v \nonumber\\
					& & (\sigma_{n} \sigma_{n-1} ..\sigma_{2} \sigma_{1} a_{n+1})^{m+1} \sigma_{n} \sigma_{n-1} .. \sigma_{n+1-r} v,\nonumber
				\end{eqnarray}
				
				where $h-1= m(n-j+1) +r $, with $0\leq r < n-j+1 $.\\
				
				So in this case (namely, $ r_{1} \leq n $), we see that  $w $ is written as
				\begin{eqnarray}
					c(\sigma_{n} \sigma_{n-1} ..\sigma_{1} a_{n+1})^{k} \sigma_{n}\sigma_{n-1} .. \sigma_{i},\nonumber
				\end{eqnarray}
				
				where  $ c$ is in $B(\tilde{A_{n-1}}) $, $ 1 \leq i \leq n+1$ and $ 0\leq k $ . \\					
		\end{itemize}

		Now, we deal with the second main case: $ r_{1} = n+1 $. Here $\sigma_{n} $ is not in the support of $\sigma_{i_{1}} ..\sigma_{2} \sigma_{1}  \sigma_{r_{1}} ..\sigma_{n-1} \sigma_{n} $ (which is equal in this case to  $\sigma_{i_{1}} ..\sigma_{2} \sigma_{1}$). Hence, we can suppose that the element in question is of the form $\sigma_{i_{1}} ..\sigma_{2} \sigma_{1}a_{n+1} w $, where $0\leq r_{0} < r_{1}$. Here we have $ r_{0} < n $, since $ r_{0} = n $ is the case of positive powers of $\sigma_{n} ..\sigma_{1} \sigma_{1}a_{n+1}$. Moreover, when $ r_{0}=0 $, then the element in question is of the form $a_{n+1} w $, which is the case $p=0$. As a consequence of this discussion we get the following corollary. \\  
	
		\begin{corollaire} \label{2_5_2}
			Let $ \bar{w} $ be fully commutative in $ W(\tilde{A_{n}}) $, where $ 2 \leq n $. Let $w$ be the corresponding element in $B(\tilde{A_{n}})$, as above. Then $w$ can be written in one and only one of the following two forms:
			\begin{eqnarray}
				& & c (\sigma_{n} \sigma_{n-1} ..\sigma_{1} a_{n+1})^{k} v, \nonumber \\ 
				\text{or }& & \sigma_{i_{0}} .. \sigma_{2}\sigma_{1}a_{n+1}c (\sigma_{n} \sigma_{n-1} ..\sigma_{1} a_{n+1})^{k} v. \nonumber	
			\end{eqnarray}
			
			Here, $c$ is in $B(\tilde{A_{n-1}})$, while $v$ is in $B(A_{n})$. \\
	
			Moreover, $ \bar{w} $ can be written in one, and only one of the following two forms (deduced from the above two forms, considering the left classes of $W(A_{n-1})$ in $W(A_{n})$): 
			\begin{eqnarray}
				& & \bar{d} (\sigma_{n} \sigma_{n-1} ..\sigma_{1} a_{n+1})^{k} \sigma_{n}\sigma_{n-1} .. \sigma_{i},\nonumber \\
				\text{or }& & \sigma_{i_{0}} .. \sigma_{2}\sigma_{1}a_{n+1} \bar{d} (\sigma_{n} \sigma_{n-1} ..\sigma_{1} a_{n+1})^{k} \sigma_{n}\sigma_{n-1} .. \sigma_{i}. \nonumber
			\end{eqnarray}
			
			Here $\bar{d}$ is in $W(\tilde{A_{n-1}})$, where  $ 1 \leq i \leq n$, with $ 0 \leq r_{0} \leq n-1 $  and $ 0\leq k $.
		\end{corollaire}
	
	
		\vspace{2cm}

\renewcommand{\refname}{REFERENCES}


\begin{thebibliography}{}

\bibitem{Graham_Lehrer_2003} J. J. Graham and G. I. Lehrer. Diagram algebras, Hecke algebras and decomposition numbers at roots of unity. Annales Scientifiques de lÉcole Normale Supérieure, 36, Issue 4:479-524, 2003. \label{Graham_Lehrer_2003}



\bibitem{Bourbaki_1981} N. Bourbaki. Groupes et algèbres de Lie, chapitres 4, 5, 6. Masson, 1981. \label{Bourbaki_1981}

\bibitem{Sadek_2013_2} S. Harbat. Markov trace on affine Temperley-Lieb algebra. 2013. \label{Sadek_2013_2}

\end{thebibliography}
\end{document}